\documentclass[hidelinks,onefignum,onetabnum]{siamart250211}

\usepackage{caption}
\usepackage{lipsum}
\usepackage{amsfonts}
\usepackage{graphicx}
\usepackage{epstopdf}
\usepackage{subcaption}
\usepackage[scr=rsfs]{mathalpha}

\usepackage{quantikz}

\usepackage{geometry}
\usepackage{amsmath,amssymb}
\usepackage{algorithm}
\usepackage{algorithmic}
\usepackage{bm}
\usepackage{appendix}
\usepackage{multirow}
\usepackage{braket}
\usepackage[english]{babel}
\usepackage{hyperref}
\usepackage[capitalize]{cleveref}
\usepackage[square,numbers,sort&compress]{natbib}
\usepackage{adjustbox}
\usepackage{xspace}
\ifpdf
  \DeclareGraphicsExtensions{.eps,.pdf,.png,.jpg}
\else
  \DeclareGraphicsExtensions{.eps}
\fi

\def\R{{\mathbb R}}
\def\N{{\mathbb N}}

\def\le{\leqslant}
\def\ge{\geqslant}

\newcommand{\eps}{\varepsilon}

\newcommand{\HI}{H_{\mathrm I}^\eps}

\numberwithin{equation}{section}

\newsiamremark{remark}{Remark}
\newsiamremark{hypothesis}{Hypothesis}
\crefname{hypothesis}{Hypothesis}{Hypotheses}
\newsiamthm{claim}{Claim}
\newsiamthm{example}{Example}

\headers{
Semiclassical Scaling Laws of Magnus Quantum Algorithms
}{Di Fang and Samuel Orellana Mateo}

\title{
Semiclassical Scaling Laws of Magnus Quantum Algorithms under Harmonic Confinement
\thanks{
\funding{This work is supported by the National Science Foundation CAREER Award DMS-2438074, and the U.S. Department of Energy, Office of Science, Accelerated Research in Quantum Computing Centers, Quantum Utility through Advanced Computational Quantum Algorithms, grant no. DE-SC0025572. }}}

\author{Di Fang\thanks{Department of Mathematics and Duke Quantum Center, Duke University 
  (\email{di.fang@duke.edu}).}
\and Samuel Orellana Mateo \thanks{Department of Mathematics, Duke University
  (\email{samuel.orellanamateo@duke.edu}).}}

\usepackage{amsopn}

\newtheorem{hyp}{Assumption}

\usepackage{url}
\usepackage{enumitem}
\usepackage{indentfirst}

\numberwithin{equation}{section}

\newcommand{\abs}[1]{\left\lvert#1\right\rvert}
\newcommand{\norm}[1]{\left\lVert#1\right\rVert}

\ifpdf
\hypersetup{
  pdftitle={Semiclassical Scaling Laws of Magnus Quantum Algorithms under Harmonic Confinement},
  pdfauthor={Di Fang and Samuel Orellana Mateo}
}
\fi

\begin{document}

\maketitle

\begin{abstract}
The simulation of semiclassical Schr\"odinger dynamics is often hindered by rapid oscillations, where classical numerical methods and quantum simulation algorithms can suffer severe precision loss as the semiclassical parameter $\eps \to 0$. In this work, we uncover a surprising robustness of high-order Magnus quantum algorithms in the semiclassical regime. We establish two important results. First, we extend the theory of Magnus superconvergence from free-particle reference dynamics to a general class of Schr\"odinger operators whose reference Hamiltonian includes harmonic confinement. Second, and more importantly, we reveal a nontrivial semiclassical scaling law: while the expansion remainder term bound suggests an error prefactor of $O(\eps^{-(p+1)})$ for a $p$th-order method, we prove that the loss is in fact order-independent and collapses to $O(\eps^{-1})$ for all $p$. The key analytical mechanism is the exact Egorov transport associated with quadratic flows in semiclassical Weyl quantization, which enables a symbol-level analysis of the interaction-picture Hamiltonian and its Magnus commutators. Our results show that high-order Magnus quantum algorithms remain far more robust in the semiclassical regime than naive series remainder term counting would suggest.
\end{abstract}



\section{Introduction}
Semiclassical Schr\"odinger equations form one of the central model classes in numerical analysis and scientific computing (see, e.g., reviews~\cite{LasserLubich2020,JinMarkowichSparber2011,BaoCai2013}, books~\cite{Lubich2008book,BlanesCasas2017book,Carles2020book}, and the papers cited herein). They arise in various applications in high-frequency quantum dynamics, molecular and chemical models, wave packet propagation, and in asymptotic descriptions that bridge quantum and classical mechanics. Their numerical treatment is challenging because the small semiclassical parameter $\eps$ produces rapid oscillations, multiple scales, and severe stability and accuracy constraints. As a result, the semiclassical regime has motivated a broad range of work on time-splitting methods, exponential integrators, geometric schemes, asymptotic-preserving ideas, and structure-aware discretizations. Quantum algorithms provide a further perspective on this problem, since they offer a natural framework for simulating oscillatory Schr\"odinger dynamics while exploiting operator structure at the Hamiltonian level.

The semiclassical regime is distinguished by pronounced multiscale structure and scale separation, and it is precisely this feature that makes both its numerical analysis and its quantum simulation subtle. Existing quantum-algorithmic studies in this regime have accordingly concentrated largely on Trotterization, or equivalently splitting methods, often together with specially structured initial states or observables~\cite{JinLiLiu2021,BornsWeilFang2022,FangTres2023,FangQu2025} for which meaningful asymptotic error estimates can still be obtained, typically inspired by the classical algorithmic line of work in the semiclassical regime, e.g.,~\cite{LasserLubich2020,JinMarkowichSparber2011,FangJinSparber2018,BaderIserlesKropielnickaSingh2014,BaderIserlesKropielnickaSingh2016,IserlesKropielnickaSingh2019,Singh2015,BaoJinMarkowich2003,BaoJakschMarkowich2003,Carles2013,CarlesGallo2017,FaouLubich2004,JinWuYang2008,FaouGradinaruLubich2009,BaoJinMarkowich2002,BaoCai2013,GolseJin19,BaoMaWang2024}. While this line of work has produced important insights, it does not exhaust the algorithmic possibilities. Other quantum simulation paradigms, especially post-Trotter approaches, are equally important from both the practical and the theoretical viewpoints. 
Among such approaches, interaction-picture methods are particularly effective when the Hamiltonian admits a decomposition into an exactly solvable dominant part and a remaining interaction. One passes to a transformed equation whose generator is time dependent, and then applies a high-order integrator adapted to that structure. The Magnus expansion has a long history in classical numerical analysis,
with significant advances on high-order integrators and their error
analysis~\cite{HochbruckLubich2003,Thalhammer2006,BlanesCasasThalhammer2017,blanes2009magnus,BlanesCasas2017book,HairerHochbruckIserlesLubich2006,ApelCubittOnorati2025}. Magnus-based quantum algorithms~\cite{FangLiuSarkar2025,FangLiuZhu2025,BornsFangZhang2025,FangZhang2025,AnFangLin2022,LiuWangZhuCongdeJongFangHuangIancuYang2026,CasaresZiniArrazola2024,SharmaTran2024,BosseChildsEtAl2024} are especially attractive in this setting: they preserve the exponential form of the propagator, encode noncommutative effects through nested commutators, and naturally interact with oscillatory dynamics. In recent years, this framework has led to a striking discovery: for certain Schr\"odinger equations in the interaction picture, a $p$th-order Magnus quantum algorithm may converge with global order $2p$, rather than merely order $p$. This doubled-order behavior, referred to as superconvergence~\cite{FangLiuSarkar2025,FangLiuZhu2025,BornsFangZhang2025,FangZhang2025,AnFangLin2022}, has now been established in several settingsand has substantially strengthened the mathematical foundation of interaction-picture quantum simulation.

The present paper is motivated by two important questions left open by this line of research. The first is an \emph{application-extension question}: can one go beyond the free-background setting and retain superconvergence when the exactly reference Hamiltonian already contains a nontrivial unbounded potential? The second is a \emph{semiclassical-dependence question}: in the equation
\begin{equation}\label{eq:semi-schd-eq-full-in-H}
i \eps \partial_t \psi^\eps(t)=H^\eps \psi^\eps(t),
\end{equation}
the interaction-picture generator carries the factor $\eps^{-1}$, and abstract local error bounds for a $p$th-order Magnus expansion involve a truncation after the first $p$ terms, whose remainder contains at least $p+1$ repetitions of the Hamiltonian. A direct counting argument therefore suggests a prefactor of size $\eps^{-(p+1)}$, which is unfavorable at high order. Is this pessimistic semiclassical scaling intrinsic, or does the operator structure of the problem enforce a sharper behavior?

To answer these questions, one must first address a nontrivial algorithmic issue. In many previous works, the natural choice of interaction-picture reference frame is the free kinetic Hamiltonian, so that the interaction term is given by the potential alone. Such a strategy is well suited to settings in which the potential is bounded. However, it is no longer appropriate in the presence of harmonic confinement, which is both physically natural and important. Indeed, for time-dependent Hamiltonian simulation, the size of the interaction-picture Hamiltonian continues to govern the algorithmic cost, regardless of whether one uses a Dyson-series method, a Magnus expansion, a discrete-clock construction, or another known simulation framework~\cite{LowWiebe2019,BerryChildsCleveEtAl2015,KieferovaSchererBerry2019,BerryChildsSuEtAl2020,AnFangLin2022,FangLiuSarkar2025,BornsFangZhang2025,FangZhang2025,FangLiuZhu2025,WatkinsWiebeRoggeroLee2022,LiWang2025,ChenGaoJiLiWangZhou2026}. Thus, if the harmonic potential were left in the interaction term, one may lose the algorithmic advantage of the interaction picture. For this reason, the first algorithmic contribution of the present work is to choose the unbounded quadratic Hamiltonian itself as the reference frame and to formulate the Magnus algorithm relative to this background.

This change of reference frame also creates new analytical difficulties. Even if one sets aside the semiclassical scaling issue and focuses only on superconvergence itself, the previous proof strategies developed for the free Schr\"odinger reference flow do not directly apply once harmonic confinement is incorporated into the reference Hamiltonian. The reason is that the interaction-picture dynamics is now generated by a quadratic background rather than by free propagation, so the associated phase-space transport and commutator structure are fundamentally different. In particular, the arguments used in the Laplacian-based setting cannot simply be transplanted to the present problem. Establishing superconvergence in this new reference frame therefore requires new numerical analysis treatment.

We then answer both questions affirmatively using this algorithmic framework. Besides the algorithmic framework, our contributions in terms of the algorithm analysis are twofold.
Our first analytical novelty is a genuine extension of the superconvergence theory itself. We prove that the doubled-order convergence phenomenon, previously known only for Laplacian-based reference dynamics, persists when the reference Hamiltonian is enlarged to include harmonic confinement. This shows that superconvergence is not a quirk of the free Schr\"odinger flow, but a structural feature of interaction-picture dynamics for a substantially richer and more physically relevant class of semiclassical problems.

Our second novelty concerns the semiclassical scaling. The transformed equation has the form
\[
\partial_t \psi^\eps_{\mathrm I}(t)
=
-\frac{i}{\eps}\, H_{\mathrm I}(t)\psi^\eps_{\mathrm I}(t).
\]
Note that here $H_{\mathrm I}$ still depends on $\eps$, but for notational convenience, we reserve the notation $H_{\mathrm I}^\eps$ for the whole Hamiltonian in the interaction picture $\HI =\eps^{-1} H_{\mathrm I}$ (see~\cref{sec:setup} for their definitions). Because the local truncation error of a $p$th-order Magnus method is governed by terms with at least $p+1$ recurrences of $H^\eps_{\mathrm I}$, one would expect the error prefactor to scale like $\eps^{-(p+1)}$. Our main result shows that this pessimistic counting is false in the present problem class. Despite the increasing commutator depth, the actual prefactor collapses sharply to order $\eps^{-1}$ for \textit{all} Magnus orders $p$, while the global convergence rate remains $2p$. Thus the algorithm error in the semiclassical scaling is \emph{order independent}. This reveals that high-order Magnus algorithms are far more robust in the semiclassical regime than abstract operator term counting alone would suggest. It also stands in contrast to other interaction-picture algorithms for time-dependent Hamiltonian simulation, such as truncated Dyson-series methods, for which a $p$th-order expansion does exhibit the expected accumulation of semiclassical factors in the algorithm error as $\eps^{-(p+1)}$.

The key analytical ingredient underlying both novelties is the exact Egorov transport associated with quadratic Hamiltonians in semiclassical Weyl quantization. In the present setting, the quadratic reference flow transports symbols exactly along the corresponding classical Hamiltonian flow. We show that the relevant symbols remain in the class $S(1)$ under the flow, which allows us to use the Calder\'on-Vaillancourt theorem to convert these symbolic bounds into operator estimates. 
This structural fact provides the appropriate analytical framework for numerical analysis. In particular, it allows us to track the interaction-picture Hamiltonian directly at the symbol level and to treat the Magnus commutators semiclassically. Through careful explicit semiclassical symbol calculations, we show that the resulting symbol remains in $S(1)$, with additional time-step convergence rate accumulation (superconvergence), yet importantly without the accumulation of powers of $\eps^{-1}$. Instead, the symbolic commutator structure collapses all orders to the same $O(\eps^{-1})$ scale. In this way, exact Egorov transport and symbolic stability explain both why superconvergence persists under harmonic confinement and why the semiclassical prefactor is far sharper than naive commutator counting predicts.

The remainder of the paper is organized as follows. In \cref{sec:setup}, we formulate the semiclassical problem, define the interaction-picture Magnus algorithm, and state the main superconvergence theorem together with its $\eps$-dependence. We prove our main results in~\cref{sec:pf_main}. In particular, \cref{sec:egorov} proves the exact Egorov theorem for the quadratic reference flow and establishes preservation of the symbol class $S(1)$. In \cref{sec:commutators}, we derive one-layer and higher nested commutator bounds and show how the naive $\eps^{-(p+1)}$ scaling collapses to $\eps^{-1}$. In \cref{sec:superconv}, we combine these estimates with the abstract Magnus error analysis to prove the main theorem. Finally, in \cref{sec:numerics}, we present numerical experiments illustrating both the doubled-order convergence and the order-independent semiclassical prefactor.

\section{Problem Setup and Main Results}\label{sec:setup}
Let $d \ge 1$ and let $x \in \mathbb{R}^d$. Let $\eps$ be the semiclassical parameter corresponding to the rescaled Planck constant.
We consider the Schr\"odinger equation in the so-called semiclassical regime, i.e., $\eps\in (0,1]$,
\begin{equation}\label{eq:schd_main}
i \eps \partial_t \psi(t,x)
=
\bigl(-\frac{\eps^2}{2}\Delta + U(x) + V(x)\bigr)\psi(t,x),
\qquad
\psi(0,x)=\psi_0(x),
\end{equation}
where the Laplacian is given by
$
\Delta = \sum_{j=1}^d \partial_{x_j}^2,
$
the potential $V$ is some smooth and bounded component, and
the potential $U : \mathbb{R}^d \to \mathbb{R}$ is a coupled harmonic oscillator
\begin{equation}
U(x) = \frac{1}{2} x^T A x = \frac{1}{2} \sum_{i,j=1}^d A_{ij} x_i x_j
\end{equation}
for a real symmetric positive definite $d \times d$ matrix $A$.
The Hamiltonian is decomposed as
\begin{equation}
H
=
\underbrace{-\frac{\eps^2}{2}\Delta + U(x)}_{H_0}
+
\underbrace{V(x)}_{V}.
\end{equation}

\begin{hyp}[Regularity of the potential $V$]\label{ass:V}
$V:\R^d\to\R$ is smooth and belongs to $S(1)$ in the sense that all derivatives of $V$ are bounded:
\[
\sup_{x\in\R^d} \abs{\partial_x^\alpha V(x)} < \infty
\qquad \text{for all multi-indices } \alpha.
\]
\end{hyp}

To design a more efficient algorithm, we pass to the interaction picture. By applying the unitary transformation $\psi_{\mathrm I}(t) = e^{i t H_0 / \eps} \psi(t)$, the unbounded reference dynamics are factored out. The evolution of the state is then governed by the time-dependent interaction Hamiltonian:
\begin{equation}\label{eq:ip_Ham_def}
    H_{\mathrm I}(t) = e^{i t H_0 / \eps} V(x) e^{-i t H_0 / \eps}.
\end{equation}
The resulting differential equation in the interaction picture is:
\begin{equation} \label{eq:interaction_picture_schd}
    i \partial_t \psi_{\mathrm I}(t) = \frac{1}{\eps} H_{\mathrm I}(t) \psi_{\mathrm I}(t).
\end{equation}
This transformation is advantageous because the interaction pictured Hamiltonian now has a smaller operator norm than the original Hamiltonian. This conversion is quantumly efficient because the reference Hamiltonian $H_0$ is exactly solvable and can be efficiently simulated quantumly (also typically called fast-forwardable; see, e.g., \cite{Somma2015,JainIyerSommaBaoJordan2025}). The full evolution is thus reduced to simulating the bounded, time-dependent operator $H_{\mathrm I}(t)$. However, notice that the generator of the dynamics in \cref{eq:interaction_picture_schd} carries a scaling factor of $1/\eps$, which is the source of the apparent semiclassical penalty.

We now focus on the interaction picture (\cref{eq:interaction_picture_schd}). The Magnus expansion provides a representation of the exact propagator $U(t) = \mathcal{T} \exp\left(-i \int_0^t \frac{1}{\eps} H_{\mathrm I}(s) \, ds\right)$ in exponential form without the time-ordering operator:
\begin{equation}
    U(t) = \exp\left( \Omega(t) \right),
\end{equation}
where the Magnus generator $\Omega(t)$ is given by the series:
\begin{align}
    \Omega(t) &= -\frac{i}{\eps} \int_0^t H_{\mathrm I}(t_1)\, dt_1 \nonumber \\
    &\quad - \frac{1}{2\eps^2} \int_0^t \int_0^{t_1} [H_{\mathrm I}(t_1), H_{\mathrm I}(t_2)] \, dt_2 dt_1 + \cdots.
\end{align}
The Magnus series converges when $\frac{1}{\eps}\int_0^t\|H(s)\|ds<\pi $ for any bounded normal operator $A(t)$. But note that we work with a finite expansion of it, and our error analysis considers the exact error representation without relying on the series expansion argument of the remainder term. Therefore, the finite sum remains well-defined. 

Based on the Magnus expansion, prior work~\cite{AnFangLin2022,FangLiuSarkar2025,FangLiuZhu2025} has established quantum Magnus algorithms that approximate the time-evolution operator over a small time step $s$ by truncating the Magnus expansion at order $p$. One uses $M$ quadrature points in each time integral, where $M$ typically needs to be as large as the spectral norm of the reference Hamiltonian $H_0$ after spatial discretization, which can be very large. Nevertheless, the quantum circuits developed in this line of work have query complexity to the Hamiltonian encoding that is independent of $M$, and elementary gate complexity that depends only polylogarithmically on $M$. Importantly, the constructions in~\cite{FangLiuZhu2025} provide efficient quantum circuits for the $p$th-order Magnus series: both the query complexity to the Hamiltonian encoding and the elementary gate complexity depend only polynomially on the order $p$, even though the high-order Magnus expansion may contain a factorial number of terms in $p$. Here, we use the same circuit construction as in~\cite{FangLiuZhu2025}, but with the interaction picture tailored to the harmonic-potential setting in the semiclassical regime. The overall gate complexity depends only polylogarithmically on the number of spatial discretization points as carefully discussed in prior work. Thus, the circuit construction itself is unchanged, and its efficiency is inherited in this case. Our focus here is instead on its error performance in the semiclassical regime, which determines how many repetitions (time-steps) of the Magnus circuit are required.

We now state our main result: the algorithmic error theorem for a $p$th-order quantum Magnus algorithm in the semiclassical regime, where $p$ is any positive integer.

\begin{theorem}[Semiclassical superconvergence with order-independent prefactor]\label{thm:main}
Let $p\in\N$ and let \cref{ass:V} hold. Let $\psi^\eps(T)$ be the exact solution to~\cref{eq:schd_main} and $\psi^{\eps,(L)}(T)$ be the $p$th-order semi-discrete Magnus approximation with time step $h=T/L$. Then there exists a constant $C_{p,T,V,A}>0$, independent of $L$ and $\eps\in(0,1]$, such that
\begin{equation}\label{eq:main-bound}
\norm{\psi^\eps(T)-\psi^{\eps,(L)}(T)}_{L^2(\R^d)}
\le C_{p,T,V,A}\, \eps^{-1} h^{2p}\, \norm{\psi_0^\eps}_{L^2(\R^d)}.
\end{equation}
\end{theorem}
In particular, the method has global order $2p$, and the semiclassical prefactor is $\eps^{-1}$ for every Magnus order $p$.
Moreover, a direct application of the Magnus remainder term counting in the celebrated Magnus convergence theory~\cite{blanes2009magnus} would suggest the weaker scaling $O(\eps^{-(p+1)}h^{p})$.
In other words, the striking point in \cref{thm:main} is not only the doubled-order convergence but also the collapse of the naive semiclassical scaling from $\eps^{-(p+1)}$ to $\eps^{-1}$. This order-independent prefactor is a direct consequence of exact Egorov transport together with semiclassical Weyl commutator structure.

\section{Proof of the Main Theorem} \label{sec:pf_main}

The main technical contributions of this section are the exact Egorov transport for the quadratic reference flow and the resulting nested commutator estimates. We organize this section as follows: in \cref{sec:preliminaries}, we briefly establish the semiclassical notation; in \cref{sec:egorov}, we establish the exact commutator lemma and the Egorov theorem for harmonic confinement; in \cref{sec:commutators}, we bound the nested commutators; and in \cref{sec:superconv}, we combine these to prove the main theorem.

\subsection{Semiclassical Notation and Preliminaries} \label{sec:preliminaries}

To establish our notation and rigorously ground the operator bounds, we briefly review the relevant tools from semiclassical analysis. 

We categorize classical observables into symbol classes. For an order function $m(x, \xi)$, the symbol class $S(m)$ consists of smooth functions $a \in C^\infty(\R^{2d})$ such that for all multi-indices $\alpha, \beta$, there exists a constant $C_{\alpha, \beta}$ satisfying $|\partial_x^\alpha \partial_\xi^\beta a(x, \xi)| \le C_{\alpha, \beta} m(x, \xi)$. As stated in \cref{ass:V}, we assume the interaction potential $V$ belongs to the standard bounded symbol class $S(1)$.

For a symbol $a \in \mathcal{S}(\R^{2d})$, its Weyl quantization is defined as the operator:
\begin{equation}
    \text{op}_\eps(a)\varphi(x) = \frac{1}{(2\pi\eps)^d} \int_{\R^{2d}} a\left(\frac{x+y}{2}, \xi\right) e^{\frac{i}{\eps}\xi \cdot (x-y)} \varphi(y) \, dy \, d\xi.
\end{equation}

Following standard semiclassical calculus \cite{zworski2022semiclassical}, the commutator of two Weyl quantized operators is exactly given by the Moyal bracket:
\begin{equation}
    [\text{op}_\eps(a), \text{op}_\eps(b)] = i\eps \, \text{op}_\eps\big( \{a, b\}_M \big),
\end{equation}
where $\{a, b\}_M(x, \xi) = \frac{2}{\eps} \sin\left( \frac{\eps}{2} P \right) \big( a(x, \xi) b(y, \eta) \big) \big|_{y=x, \eta=\xi}$, and $P = \langle \partial_x, \partial_\eta \rangle - \langle \partial_\xi, \partial_y \rangle$ is the symplectic bidifferential operator.

Finally, to bound the operator norms uniformly in the semiclassical parameter, we rely on the Calder\'{o}n-Vaillancourt theorem \cite{zworski2022semiclassical, Martinez2002}. It guarantees that if $a \in S(1)$, its Weyl quantization extends to a bounded operator on $L^2(\R^d)$ satisfying:
\begin{equation}
    \|\text{op}_\eps(a)\|_{L^2 \to L^2} \le C_d \sum_{|\alpha| + |\beta| \le M_d} \sup_{(x, \xi) \in \R^{2d}} \left| \partial_x^\alpha \partial_\xi^\beta a(x, \xi) \right|,
\end{equation}
where the constants $C_d$ and $M_d$ depend only on the dimension $d$ and are strictly independent of $\eps \in (0, 1]$.

\subsection{Exact Egorov Transport for Harmonic Confinement} \label{sec:egorov}

To analyze the algorithmic error of the Magnus expansion in the semiclassical regime, we must first understand the structure of the interaction picture Hamiltonian $H_{\mathrm{I}}(t) = e^{itH_0/\eps} V(x) e^{-itH_0/\eps}$. 

We begin by explicitly solving the classical flow associated with the reference Hamiltonian. Let the classical phase space be $\R^{2d} = \R^d_x \times \R^d_\xi$. The classical reference Hamiltonian corresponding to $H_0$ is the quadratic form:
\begin{equation}
    h_0(x, \xi) = \frac{1}{2}|\xi|^2 + \frac{1}{2}x^T A x.
\end{equation}
The classical flow $\Phi_t : \R^{2d} \to \R^{2d}$ generated by $h_0$ is determined by Hamilton's equations:
\begin{equation}
    \frac{d}{dt} \begin{pmatrix} x \\ \xi \end{pmatrix} = \begin{pmatrix} 0 & I \\ -A & 0 \end{pmatrix} \begin{pmatrix} x \\ \xi \end{pmatrix} = M \begin{pmatrix} x \\ \xi \end{pmatrix}.
\end{equation}
The exact solution is given by the matrix exponential $\Phi_t(x, \xi) = \exp(tM) (x, \xi)^T$. By evaluating the powers of the block matrix $M$, we observe the pattern $M^{2m} = (-1)^m \text{diag}(A^m, A^m)$ and $M^{2m+1} = (-1)^m \begin{pmatrix} 0 & A^m \\ -A^{m+1} & 0 \end{pmatrix}$. Because $A$ is symmetric and positive-definite, it admits a unique symmetric positive-definite square root matrix $\Omega = A^{1/2}$, such that $A^m = \Omega^{2m}$. Substituting these into the Taylor series expansion of the matrix exponential yields the explicit classical flow operator:
\begin{equation}
    e^{tM} = \begin{pmatrix} \cos(\Omega t) & \Omega^{-1}\sin(\Omega t) \\ -\Omega\sin(\Omega t) & \cos(\Omega t) \end{pmatrix}.
\end{equation}
Thus, the classical trajectories $\Phi_t(x, \xi) = (x(t), \xi(t))$ are given exactly by:
\begin{equation} \label{eq:classical_flow_explicit}
    \begin{cases}
        x(t) = \cos(\Omega t)x + \Omega^{-1}\sin(\Omega t)\xi \\
        \xi(t) = -\Omega\sin(\Omega t)x + \cos(\Omega t)\xi.
    \end{cases}
\end{equation}
Notice that $\Phi_t$ is a linear transformation with bounded oscillatory coefficients. Furthermore, because it is generated by a Hamiltonian, $\Phi_t$ is a canonical symplectic mapping for all $t \in \R$.

Next, we establish a fundamental algebraic identity linking the quantum commutator of the reference Hamiltonian to the classical Poisson bracket.

\begin{lemma}[Exact Commutator] \label{lem:exact_commutator}
Let $H_0 = \text{op}_\eps(h_0) = -\frac{\eps^2}{2}\Delta_x + \frac{1}{2}x^T A x$ be the Weyl quantization of the harmonic oscillator Hamiltonian, and let $a \in \mathcal{S}(\R^{2d})$ be a Schwartz class symbol. Then the commutator is exactly given by the Weyl quantization of the Poisson bracket:
\begin{equation}
    [H_0, \text{op}_\eps(a)] = \frac{\eps}{i} \text{op}_\eps(\{h_0, a\}),
\end{equation}
where $\{h_0, a\} = \nabla_\xi h_0 \cdot \nabla_x a - \nabla_x h_0 \cdot \nabla_\xi a$.
\end{lemma}

\begin{proof}[Proof of \cref{lem:exact_commutator}]
We establish the identity for a Schwartz symbol $a\in\mathcal S(\mathbb R^{2d})$ and a test function $\varphi\in\mathcal S(\mathbb R^d)$. We decompose the Hamiltonian into its kinetic and potential components, writing $H_0=K+U$, where $K=-\frac{\eps^2}{2}\Delta_x$ and $U(x)=\frac12 x^T A x$.

We first compute the kinetic contribution. By definition, the action of $\operatorname{op}_\eps(a)$ on the test function is given by the integral
\[
\operatorname{op}_\eps(a)\varphi(x) 
= 
\frac{1}{(2\pi\eps)^d} \int_{\mathbb R^{2d}} a(q,\xi) E(x,y,\xi) \varphi(y) \, dy \, d\xi,
\]
where we define the midpoint $q=(x+y)/2$ and the oscillatory exponential $E(x,y,\xi) := \exp\!\left(\frac{\mathrm i}{\eps}\xi\cdot(x-y)\right)$. To evaluate the commutator with the Laplacian, we apply $\Delta_x$ to the integrand $a(q,\xi)E$. Using the product rule $\Delta(uv) = (\Delta u)v + 2\nabla u \cdot \nabla v + u(\Delta v)$ along with the chain rule relations $\nabla_x a = \frac{1}{2}\nabla_q a$ and $\Delta_x a = \frac{1}{4}\Delta_q a$, we obtain
\[
\Delta_x\bigl(a(q,\xi)E\bigr) 
= 
\frac{1}{4}(\Delta_q a)E + \frac{2\mathrm i}{\eps}\left(\frac{1}{2}\nabla_q a \cdot \xi\right)E - a\frac{|\xi|^2}{\eps^2}E.
\]
Similarly, applying the Laplacian with respect to $y$, and noting that $\nabla_y a = \frac{1}{2}\nabla_q a$ but $\nabla_y E = -\frac{\mathrm i}{\eps}\xi E$, yields
\[
\Delta_y\bigl(a(q,\xi)E\bigr) 
= 
\frac{1}{4}(\Delta_q a)E - \frac{2\mathrm i}{\eps}\left(\frac{1}{2}\nabla_q a \cdot \xi\right)E - a\frac{|\xi|^2}{\eps^2}E.
\]
Subtracting the two expressions causes the symmetric second-order and zero-order terms to cancel exactly, leaving 
\[
\Delta_x\bigl(a(q,\xi)E\bigr) - \Delta_y\bigl(a(q,\xi)E\bigr)
=
\frac{2\mathrm i}{\eps}(\nabla_q a(q,\xi)\cdot\xi)E.
\]
Substituting this identity back into the integral definition of $\Delta_x\operatorname{op}_\eps(a)\varphi(x)$, we have
\[
\Delta_x \operatorname{op}_\eps(a)\varphi(x) 
= 
\frac{1}{(2\pi\eps)^d} \int_{\mathbb R^{2d}} \left( \Delta_y \bigl(a(q, \xi) E\bigr) + \frac{2\mathrm i}{\eps} (\nabla_q a \cdot \xi) E \right) \varphi(y) \, dy \, d\xi.
\]
Because $\varphi$ and $a$ exhibit rapid Schwartz decay, we can integrate the first term by parts with respect to $y$. All boundary terms vanish, shifting the Laplacian onto $\varphi(y)$ and yielding $\int \Delta_y(aE)\varphi \, dy = \int aE(\Delta_y \varphi) \, dy$. Reconstructing the operators, this gives
\[
\Delta_x\operatorname{op}_\eps(a)\varphi
=
\operatorname{op}_\eps(a)\Delta_x\varphi
+
\frac{2\mathrm i}{\eps}\operatorname{op}_\eps(\xi\cdot\nabla_x a)\varphi.
\]
Multiplying by $-\frac{\eps^2}{2}$ and rearranging terms isolates the kinetic commutator:
\[
\left[-\frac{\eps^2}{2}\Delta_x,\operatorname{op}_\eps(a)\right]
=
\frac{\eps}{\mathrm i}\operatorname{op}_\eps(\xi\cdot\nabla_x a).
\]

Next, we compute the potential contribution by evaluating the difference $[U,\operatorname{op}_\eps(a)]\varphi(x) = U(x)\operatorname{op}_\eps(a)\varphi(x) - \operatorname{op}_\eps(a)(U\varphi)(x)$. Bringing the potential directly into the integral, we get
\[
[U, \operatorname{op}_\eps(a)]\varphi(x) 
= 
\frac{1}{(2\pi\eps)^d} \int_{\mathbb R^{2d}} a(q, \xi) \bigl(U(x) - U(y)\bigr) E(x,y,\xi) \varphi(y) \, dy \, d\xi.
\]
To evaluate $U(x)-U(y)$, we express $x$ and $y$ in terms of the midpoint $q$ and the difference $\delta=x-y$, so that $x=q+\delta/2$ and $y=q-\delta/2$. Expanding the quadratic form $U(x) = \frac{1}{2}x^T A x$, we have
\[
U(x) = \frac{1}{2}\left(q + \frac{\delta}{2}\right)^T A \left(q + \frac{\delta}{2}\right) = \frac{1}{2}q^T A q + \frac{1}{2}q^T A \delta + \frac{1}{8}\delta^T A \delta,
\]
and similarly $U(y) = \frac{1}{2}q^T A q - \frac{1}{2}q^T A \delta + \frac{1}{8}\delta^T A \delta$. Subtracting the two, the quadratic $\mathcal{O}(\delta^2)$ terms cancel perfectly with no higher-order remainder, yielding the exact identity $U(x)-U(y) = q^T A \delta = \nabla_x U(q)\cdot(x-y)$. Substituting this exact difference back into the integral gives
\[
[U, \operatorname{op}_\eps(a)]\varphi(x) 
= 
\frac{1}{(2\pi\eps)^d} \int_{\mathbb R^{2d}} a(q, \xi) \bigl(\nabla_x U(q) \cdot (x-y)\bigr) E(x,y,\xi) \varphi(y) \, dy \, d\xi.
\]
We now utilize the algebraic relation $(x-y)E = \frac{\eps}{\mathrm i}\nabla_\xi E$. Integration by parts with respect to $\xi$ (where boundary terms again vanish due to the Schwartz decay of $a$) shifts the momentum gradient onto $a(q,\xi)$, pulling out a minus sign and resulting in
\[
[U,\operatorname{op}_\eps(a)]
=
-\frac{\eps}{\mathrm i}
\operatorname{op}_\eps(\nabla_\xi a\cdot\nabla_xU).
\]

Combining the kinetic and potential commutator identities, we obtain the full commutator:
\[
[H_0,\operatorname{op}_\eps(a)]
=
\frac{\eps}{\mathrm i}
\operatorname{op}_\eps
\bigl(
\xi\cdot\nabla_x a
-
\nabla_xU\cdot\nabla_\xi a
\bigr).
\]
Because the classical Hamiltonian is $h_0(x,\xi) = \frac{1}{2}|\xi|^2 + U(x)$, its phase space gradients are precisely $\nabla_\xi h_0=\xi$ and $\nabla_x h_0=\nabla_xU$. The expression inside the quantization is therefore exactly the classical Poisson bracket $\{h_0,a\} = \nabla_\xi h_0 \cdot \nabla_x a - \nabla_x h_0 \cdot \nabla_\xi a$, which completes the proof.
\end{proof}

With the exact commutator identity established, we can now prove that the quantum evolution generated by the harmonic oscillator exactly transports symbols along the classical flow $\Phi_t$.

\begin{theorem}[Exact Egorov Theorem for Quadratic Hamiltonians] \label{thm:exact_egorov}
Let $H_0 = \text{op}_\eps(h_0)$ be the Weyl quantized operator of the harmonic oscillator $h_0(x, \xi) = \frac{1}{2}|\xi|^2 + \frac{1}{2}x^T A x$, and let $\Phi_t : \R^{2d} \to \R^{2d}$ be the classical flow generated by $h_0$. For any symbol $a \in S(1)$ and every $t\in\R$, the Heisenberg evolution of the quantized observable is exactly the quantization of the classically evolved observable:
\begin{equation}
    e^{itH_0/\eps} \text{op}_\eps(a) e^{-itH_0/\eps} = \text{op}_\eps(a \circ \Phi_t).
\end{equation}
\end{theorem}

\begin{proof}[Proof of \cref{thm:exact_egorov}]
We first establish the identity for a Schwartz symbol $a\in\mathcal S(\mathbb R^{2d})$. We then extend the result to $S(1)$ by a cutoff argument.

Assuming first that $t\ge0$ (the case $t<0$ is identical under time reversal), we track the quantum evolution by defining the time-dependent interpolating operator
\[
A(s)
:=
e^{\mathrm i sH_0/\eps}
\operatorname{op}_\eps(a\circ\Phi_{t-s})
e^{-\mathrm i sH_0/\eps},
\qquad
0\le s\le t.
\]
By construction, evaluating at the endpoints gives the target evolved operator $A(t) = e^{\mathrm i tH_0/\eps} \operatorname{op}_\eps(a) e^{-\mathrm i tH_0/\eps}$ and the classically transported observable $A(0) = \operatorname{op}_\eps(a\circ\Phi_t)$. To establish that these operators are identical, we invoke the Fundamental Theorem of Calculus, $A(t) - A(0) = \int_0^t A'(s) \, ds$, and will show that the derivative $A'(s)=0$ for all $s\in[0,t]$.

Applying the product rule to the three time-dependent factors defining $A(s)$, we obtain
\[
\begin{aligned}
A'(s)
&=
\left(\frac{\mathrm i}{\eps}H_0 e^{\mathrm i sH_0/\eps}\right)
\operatorname{op}_\eps(a\circ\Phi_{t-s})
e^{-\mathrm i sH_0/\eps}
+
e^{\mathrm i sH_0/\eps}
\operatorname{op}_\eps\bigl(\partial_s(a\circ\Phi_{t-s})\bigr)
e^{-\mathrm i sH_0/\eps}
\\
&\quad+
e^{\mathrm i sH_0/\eps}
\operatorname{op}_\eps(a\circ\Phi_{t-s})
\left(-\frac{\mathrm i}{\eps}H_0 e^{-\mathrm i sH_0/\eps}\right).
\end{aligned}
\]
By factoring out the unitary propagators $e^{\mathrm i sH_0/\eps}$ on the left and $e^{-\mathrm i sH_0/\eps}$ on the right, the remaining terms natively group into a quantum commutator:
\[
A'(s)
=
e^{\mathrm i sH_0/\eps}
\Bigl(
\frac{\mathrm i}{\eps}
[H_0,\operatorname{op}_\eps(a\circ\Phi_{t-s})]
+
\operatorname{op}_\eps\bigl(\partial_s(a\circ\Phi_{t-s})\bigr)
\Bigr)
e^{-\mathrm i sH_0/\eps}.
\]
According to \cref{lem:exact_commutator}, the commutator term evaluates exactly to the Weyl quantization of the Poisson bracket, yielding $\operatorname{op}_\eps\bigl(\{h_0,a\circ\Phi_{t-s}\}\bigr)$. 

It remains to explicitly evaluate the $s$-derivative of the transported symbol. Let $\tau=t-s$. By the multivariate chain rule, the total time derivative of the observable $a$ evaluated along the classical trajectory $\Phi_\tau(x,\xi)$ is given by
\[
\partial_\tau(a\circ\Phi_\tau)
=
\nabla_x a(\Phi_\tau)\cdot \dot{x}(\tau)
+
\nabla_\xi a(\Phi_\tau)\cdot \dot{\xi}(\tau).
\]
Substituting Hamilton's equations of motion, $\dot{x} = \nabla_\xi h_0$ and $\dot{\xi} = -\nabla_x h_0$, directly into this expression, we find
\[
\partial_\tau(a\circ\Phi_\tau)
=
\nabla_x a(\Phi_\tau)\cdot\nabla_\xi h_0(\Phi_\tau)
-
\nabla_\xi a(\Phi_\tau)\cdot\nabla_x h_0(\Phi_\tau)
=
\{h_0,a\}\circ\Phi_\tau.
\]
Because $\Phi_\tau$ is a canonical symplectic transformation generated by $h_0$, the Hamiltonian $h_0$ is a conserved quantity, and the Poisson bracket is invariant under the flow. This invariance implies $\{h_0,a\}\circ\Phi_\tau = \{h_0,a\circ\Phi_\tau\}$. Applying the chain rule for $s$, we consequently find
\[
\partial_s(a\circ\Phi_{t-s})
=
-\partial_\tau(a\circ\Phi_{\tau}) \big|_{\tau=t-s}
=
-\{h_0,a\circ\Phi_{t-s}\}.
\]
Substituting this relation back into the expression for $A'(s)$, the Poisson bracket generated by the commutator and the derivative of the classical flow cancel exactly inside the quantization. Thus, $A'(s)=0$ for all $s\in[0,t]$. Because the derivative is identically zero, the integral from the Fundamental Theorem of Calculus vanishes, implying $A(t)=A(0)$. This proves the identity for Schwartz symbols.

To extend the result to $a\in S(1)$, choose $\chi\in C_c^\infty(\mathbb R^{2d})$ equal to one near the origin and define
\[
a_R(z)=\chi(z/R)a(z), \qquad R\ge1.
\]
Then $a_R\in\mathcal S(\mathbb R^{2d})$ and $a_R\to a$ pointwise, while the family $(a_R)$ is uniformly bounded in $S(1)$. Since $\Phi_t$ is linear, the same properties hold for $a_R\circ\Phi_t$ and $a\circ\Phi_t$.

For $f,g\in\mathcal S(\mathbb R^d)$, the Weyl--Wigner pairing gives
\[
\langle\operatorname{op}_\eps(a_R)f,g\rangle
=
\int_{\mathbb R^{2d}}
a_R(z)W_\eps(f,g)(z)\,dz,
\]
where $W_\eps(f,g)$ is the cross-Wigner transform. Because $W_\eps(f,g)$ is Schwartz and the symbols $a_R$ are uniformly bounded, dominated convergence implies
\[
\langle\operatorname{op}_\eps(a_R)f,g\rangle
\longrightarrow
\langle\operatorname{op}_\eps(a)f,g\rangle.
\]
The same argument applies to the transported symbols $a_R\circ\Phi_t$.

The harmonic-oscillator propagator preserves $\mathcal S(\mathbb R^d)$. We may therefore pass to the limit in
\[
e^{itH_0/\eps}\operatorname{op}_\eps(a_R)e^{-itH_0/\eps}
=
\operatorname{op}_\eps(a_R\circ\Phi_t)
\]
after pairing both sides with arbitrary Schwartz functions $f,g$. Finally, the Calder\'on-Vaillancourt theorem ensures that both limiting operators are bounded on $L^2(\mathbb R^d)$, so the identity extends by density of Schwartz test functions in $L^2$ to all $L^2(\mathbb R^d)$.

\end{proof}

\subsection{Nested Commutator Estimates} \label{sec:commutators}

By combining the exact Egorov transport (\cref{thm:exact_egorov}) with the semiclassical Moyal bracket, we can track the interaction-picture Hamiltonian directly at the symbol level. Because the classical flow $\Phi_t$ is a linear transformation with bounded coefficients, the transported symbol $V \circ \Phi_t$ remains in the class $S(1)$. This symbolic stability, together with the Calder\'{o}n-Vaillancourt theorem, allows us to rigorously bound the operator norms of the nested commutators and reveal the collapse of the semiclassical scaling.

\begin{theorem}[Nested Commutator Estimate]
\label{thm:nested_comm}
Let $V \in S(1)$ satisfy \cref{ass:V}, and let $H_{\mathrm{I}}(t) = e^{itH_0/\eps} V(x) e^{-itH_0/\eps}$ be the interaction-picture Hamiltonian. Let $h > 0$ be the time step size, and consider a time interval $[t_j, t_{j+1}]$ of length $h$. For any integer $q \ge 2$, let $\mathscr{C}_q$ denote the set of all nested commutators of grade $q$ (containing exactly $q-1$ commutator brackets). Then, there exists a constant $C_{q,T} > 0$, depending only on $q$, the final time $T$, the matrix $A$, and the $S(1)$ seminorms of $V$, such that
\begin{equation}
    \sup_{\tau_1, \dots, \tau_q \in[t_j, t_{j+1}]} \max_{C \in \mathscr{C}_q} \left\| C\left( \frac{1}{\eps}H_{\mathrm{I}}(\tau_1), \dots, \frac{1}{\eps}H_{\mathrm{I}}(\tau_q) \right) \right\|_{L^2 \to L^2} \le C_{q,T} \, \eps^{-1} h^{q-1}.
\end{equation}
\end{theorem}

\begin{proof}[Proof of \cref{thm:nested_comm}]
We define the classically transported symbol $v_t:=V\circ\Phi_t$. By \cref{thm:exact_egorov}, the interaction-picture Hamiltonian has the exact Weyl representation $H_{\mathrm I}(t)=\operatorname{op}_\eps(v_t)$. Since $\Phi_t$ is linear with uniformly bounded coefficients, the symbols $v_t$ are uniformly bounded in $S(1)$.

We fix the base time $t_*=t_j$ and write $s_i=\tau_i-t_*\in[0,h]$. The identity
\[
H_{\mathrm I}(t_*+s)=e^{i  t_*H_0/\eps}H_{\mathrm I}(s)e^{-i  t_*H_0/\eps}
\]
and invariance of the operator norm under unitary conjugation reduce the estimate to times in $[0,h]$.

First suppose $0<h\le1$. In what follows, we consider $\eps h$ as the effective semiclassical parameter.
Rescaling the momentum variable in Weyl quantization gives
\[
H_{\mathrm I}(s)=\operatorname{op}_{\eps h}(v_{s,h}),\qquad
v_{s,h}(x,\xi):=v_s(x,\xi/h)
=V\!\left(\cos(\Omega s)x+
h^{-1}\Omega^{-1}\sin(\Omega s)\xi\right).
\]

For $0\le s\le h$, we have
\[
\norm{\cos(\Omega s)}\le1,\qquad
\norm{h^{-1}\Omega^{-1}\sin(\Omega s)}\le s/h\le1.
\]
It follows that all $S(1)$ seminorms of $v_{s,h}$ are bounded uniformly in $s$, $h$, and $\eps$.

To evaluate the nested commutators, we use the Moyal-bracket identity with semiclassical parameter $\eps h$,
\[
[\operatorname{op}_{\eps h}(a),\operatorname{op}_{\eps h}(b)]
=i \eps h\,\operatorname{op}_{\eps h}(\{a,b\}_M).
\]
The $S(1)$ seminorms of this Moyal bracket are controlled by finitely many seminorms of $a$ and $b$, uniformly for $0<\eps h\le1$. Repeated application of this identity therefore gives, for any $C\in\mathscr C_q$,
\[
C\bigl(\eps^{-1}H_{\mathrm I}(s_1),\ldots,
\eps^{-1}H_{\mathrm I}(s_q)\bigr)
=\eps^{-q}(\eps h)^{q-1}\operatorname{op}_{\eps h}(c_q),
\]
where $c_q$ is uniformly bounded in $S(1)$. This applies to every bracket structure: if the two subcommutators contain $r$ and $q-r$ entries, their factors $(\eps h)^{r-1}$ and $(\eps h)^{q-r-1}$ acquire one further factor $\eps h$ from the outer bracket.

Applying the Calder\'on--Vaillancourt theorem, we obtain
\[
\norm{C\bigl(\eps^{-1}H_{\mathrm I}(s_1),\ldots,
\eps^{-1}H_{\mathrm I}(s_q)\bigr)}_{L^2\to L^2}
\le C_{q,T}\eps^{-q}(\eps h)^{q-1}
=C_{q,T}\eps^{-1}h^{q-1}.
\]
If $h>1$, the same argument with parameter $\eps$ and the uniformly bounded symbols $v_s$ gives $C_{q,T}\eps^{-1}\le C_{q,T}\eps^{-1}h^{q-1}$. Taking the supremum over the times and the maximum over the finitely many grade-$q$ bracket structures completes the proof.
\end{proof}

\subsection{Superconvergence and Global Error Bound} \label{sec:superconv}

Combining the nested commutator estimates with the general local error theory for Magnus expansions and the standard accumulation of local errors under unitary evolution immediately gives~\cref{thm:main}.

\begin{lemma}[Local truncation error for the $p$th-order Magnus method~\cite{FangLiuZhu2025}]\label{lem:local-magnus}
Let $U_\mathrm{exact}(t_{j+1},t_j)$ be the exact propagator and let $U_p(t_{j+1},t_j)$ be the $p$th-order Magnus approximation over a step of length $h = T/L$ and $t_j = j h$. Then
\begin{equation}\label{eq:local-template}
\norm{U_\mathrm{exact}(t_{j+1},t_j)-U_p(t_{j+1},t_j)}_{L^2\to L^2}
\le C_p \sum_{q=p+1}^{p^2+2p} \alpha_{\mathrm{comm},q}\, h^q,
\end{equation}
where $\alpha_{\mathrm{comm},q}$ is the maximal norm of grade-$q$ nested commutators (i.e. $\mathscr C_q$) of $\varepsilon^{-1} H_\mathrm{I}(t)$ on $[t_j,t_{j+1}]$.
\end{lemma}
The estimate applies to the present bounded, strongly continuous family. Indeed, norm-continuous convex interpolants of finitely many time samples preserve the commutator bounds by multilinearity. Their exact propagators and finite Magnus approximations converge strongly, so the same operator-norm bound follows by passing to the limit on each vector.

\begin{proof}[Proof of \cref{thm:main}]
By \cref{thm:nested_comm}, the commutator factor in \cref{eq:local-template} satisfies
\begin{equation}
    \alpha_{\mathrm{comm},q} \le C_{q,T} \, \eps^{-1} h^{q-1}.
\end{equation}
Substituting this into \cref{eq:local-template}, we bound the local truncation error. Expanding the sum, we have:
\begin{align}
    \left\| U_{\mathrm{exact}}(t_{j+1}, t_j) - U_p(t_{j+1}, t_j) \right\|_{L^2 \to L^2} 
    &\le C_p \sum_{q=p+1}^{p^2+2p} \left( C_{q,T} \, \eps^{-1} h^{q-1} \right) h^q \nonumber \\
    &= \eps^{-1} \sum_{q=p+1}^{p^2+2p} \tilde{C}_{q,T} \, h^{2q-1},
\end{align}
where $\tilde{C}_{q,T} = C_p C_{q,T}$. For $h \le 1$, the polynomial in $h$ is strictly dominated by its lowest-degree term, which is $h^{2p+1}$ (corresponding to $q=p+1$). We can therefore factor out $h^{2p+1}$ to obtain a rigorous local error bound:
\begin{equation}
    \left\| U_{\mathrm{exact}}(t_{j+1}, t_j) - U_p(t_{j+1}, t_j) \right\|_{L^2 \to L^2} \le C_{\mathrm{local}} \, \eps^{-1} h^{2p+1},
\end{equation}
where $C_{\mathrm{local}}$ depends only on $p$, $T$, $V$, and $A$.

To transition from the local error over a single step to the global error over the entire simulation interval $[0, T]$, we use a telescoping sum identity. Let $U_{\mathrm{num}}(t_k, t_j)$ denote the numerical evolution obtained by composing the single-step Magnus propagators. We define a sequence of hybrid operators $E_j$ for $j = 0, 1, \dots, L$:
\begin{equation}
    E_j = U_{\mathrm{num}}(T, t_j) U_{\mathrm{exact}}(t_j, 0).
\end{equation}
Notice the boundary cases $E_L = U_{\mathrm{exact}}(T, 0)$ and $E_0 = U_{\mathrm{num}}(T, 0)$. Therefore, the global error operator can be written exactly as the sum of the differences of consecutive hybrid operators:
\begin{equation}
    U_{\mathrm{exact}}(T, 0) - U_{\mathrm{num}}(T, 0) = E_L - E_0 = \sum_{j=0}^{L-1} (E_{j+1} - E_j).
\end{equation}
We expand the difference $E_{j+1} - E_j$:
\begin{align}
    E_{j+1} - E_j &= U_{\mathrm{num}}(T, t_{j+1}) U_{\mathrm{exact}}(t_{j+1}, 0) - U_{\mathrm{num}}(T, t_j) U_{\mathrm{exact}}(t_j, 0) \nonumber \\
    &= U_{\mathrm{num}}(T, t_{j+1}) \Big[ U_{\mathrm{exact}}(t_{j+1}, t_j) - U_p(t_{j+1}, t_j) \Big] U_{\mathrm{exact}}(t_j, 0).
\end{align}
Taking the $L^2 \to L^2$ operator norm of both sides, and utilizing the submultiplicativity of the norm along with the exact unitarity of the propagators (i.e., $\|U_{\mathrm{num}}(T, t_{j+1})\| = 1$ and $\|U_{\mathrm{exact}}(t_j, 0)\| = 1$), we obtain:
\begin{equation}
    \left\| E_{j+1} - E_j \right\|_{L^2 \to L^2} \le 1 \cdot \left\| U_{\mathrm{exact}}(t_{j+1}, t_j) - U_p(t_{j+1}, t_j) \right\|_{L^2 \to L^2} \cdot 1 \le C_{\mathrm{local}} \, \eps^{-1} h^{2p+1}.
\end{equation}
Applying the triangle inequality to the telescoping sum, the global operator error is bounded by the sum of the local errors:
\begin{align}
    \left\| U_{\mathrm{exact}}(T, 0) - U_{\mathrm{num}}(T, 0) \right\|_{L^2 \to L^2} 
    &\le \sum_{j=0}^{L-1} \left\| E_{j+1} - E_j \right\|_{L^2 \to L^2} \nonumber \\
    &\le L \left( C_{\mathrm{local}} \, \eps^{-1} h^{2p+1} \right) = \left( \frac{T}{h} \right) C_{\mathrm{local}} \, \eps^{-1} h^{2p+1} = C_{p,T,V,A} \, \eps^{-1} h^{2p},
\end{align}
where $C_{p,T,V,A} = C_{\mathrm{local}} T$. 

Finally, we apply this operator bound to the initial quantum state $\psi_0^\eps$. Let $\psi_{\mathrm{I}}^\eps(T) = U_{\mathrm{exact}}(T, 0) \psi_0^\eps$ be the exact state in the interaction picture at time $T$, and let $\psi_{\mathrm{I}}^{\eps,(L)}(T) = U_{\mathrm{num}}(T, 0) \psi_0^\eps$ be the numerical state. By the definition of the operator norm, we have:
\begin{equation}
    \left\| \psi_{\mathrm{I}}^\eps(T) - \psi_{\mathrm{I}}^{\eps,(L)}(T) \right\|_{L^2(\R^d)} \le C_{p,T,V,A} \, \eps^{-1} h^{2p} \left\| \psi_0^\eps \right\|_{L^2(\R^d)}.
\end{equation}
Because the transformation back to the Schr\"{o}dinger picture involves multiplying by the unitary operator $e^{-i T H_0 / \eps}$, which preserves the $L^2$ norm exactly, the error bound holds identically for the full Schr\"{o}dinger wavefunctions $\psi^\eps(T)$ and $\psi^{\eps,(L)}(T)$.

\end{proof}

\section{Circuit Description}\label{sec:circuit}
Here, we use the same interaction-picture circuit framework as in~\cite{AnFangLin2022,FangLiuSarkar2025,LowWiebe2019}, together with the Magnus circuit construction in~\cite{FangLiuSarkar2025,FangLiuZhu2025}, but with the reference evolution tailored to the harmonic-potential setting in the semiclassical regime. For the circuit discussion, we restrict attention to the one-dimensional harmonic oscillator
\[
H_0 =-\frac{\eps^2}{2}\Delta+\frac12 \omega^2 x^2, \quad \omega>0.
\]
Efficient circuit constructions for quantum harmonic oscillators without semiclassical scaling have already been studied in the literature. Our purpose here is not to re-develop those general constructions, but rather to show, in this simple and representative setting, how the interaction-picture circuit can be adapted to the semiclassical regime. For more general quantum harmonic oscillators and their fast-forwarding constructions, we refer the reader to~\cite{JainIyerSommaBaoJordan2025}, which can also be adapted to the semiclassical regime following the same procedure described here. Since the emphasis of the present paper is on the numerical-analysis side, we do not repeat those broader circuit details here.

A key observation is that the semiclassical harmonic-oscillator propagator admits an explicit factorization. For $|\omega t|<\pi$,
\[
e^{-itH_0/\eps}
=
e^{-ia(t)x^2/\eps}
e^{-ib(t)p_\eps^2/\eps}
e^{-ia(t)x^2/\eps},
\qquad
p_\eps=-i\eps\partial_x,
\]
where
\[
a(t)=\frac{\omega}{2}
\tan\!\left(\frac{\omega t}{2}\right),
\qquad
b(t)=\frac{1}{2\omega}\sin(\omega t).
\]

Write $U_0(t)=e^{-i tH_0/\eps}$. The coefficients $a(t)$ and $b(t)$ are $2\pi/\omega$-periodic,
whereas $U_0(t+2\pi/\omega)=-U_0(t)$. Hence, whenever
$\cos(\omega t/2)\ne0$, the three-factor product agrees with
$U_0(t)$ up to a scalar sign, even outside the short-time interval. The interaction-picture Hamiltonian is $U_0(t)^\dagger VU_0(t)$. Therefore, any scalar phase in the reference evolution $U_0(t)$ does not affect the result: if $\widetilde U_0(t)=e^{i\theta(t)}U_0(t)$, then
\[
\widetilde U_0(t)^\dagger V\widetilde U_0(t)
=U_0(t)^\dagger VU_0(t)=H_{\mathrm I}(t).
\]
The same cancellation holds for a coherent superposition of time indices, since the time-dependent phase acts only on the time register and commutes with the potential block. The inverse reference evolution is implemented as the adjoint of the same circuit.

Since $U_0(t+2\pi/\omega)=-U_0(t)$, we may reduce each queried time modulo $2\pi/\omega$ to $t'\in[-\pi/\omega,\pi/\omega)$. We then use
\[
U_0(t')=\left(e^{-i(t'/2)H_0/\eps}\right)^2,
\qquad |\omega t'/2|\le\pi/2.
\]
Thus two short-time blocks suffice, with the tangent arguments uniformly away from their poles. For coherently queried times $t_m$, the reduction to $t'_m$ and the coefficient evaluation are performed coherently. The number of blocks is independent of $t$, $\eps$, and $N$; the arithmetic cost depends on the bit length and required precision of the time input.

We now describe one short-time block, assuming $|\omega t|<\pi$ and, for its controlled version, $|\omega t_m|<\pi$. For an arbitrary queried time, the interaction-picture construction uses two such blocks at $t'_m/2$.

\begin{figure}[t]
    \centering
    \begin{subfigure}[t]{0.95\textwidth}
        \centering
\[
\begin{quantikz}[row sep=0.35cm, column sep=0.65cm]
\lstick{$\ket{m}_d$}
& \gate[2]{\mathcal A_{ab}}
& \qw
& \qw
& \qw
& \qw
& \qw
& \gate[2]{\mathcal A_{ab}^{\dagger}}
& \qw
\\
\lstick{$\ket{0}_c$}
&
& \gate[2]{\mathrm{Ph}_x}
& \qw
& \gate[2]{\mathrm{Ph}_{p_\eps}}
& \qw
& \gate[2]{\mathrm{Ph}_x}
&
& \qw
\\
\lstick{$\ket{\psi}_s$}
& \qw
&
& \gate{F_N}
&
& \gate{F_N^\dagger}
&
& \qw
\end{quantikz}
\]
        
     \caption{Circuit for the reference propagator $U_{0,d}(t_m)\approx e^{-it_mH_0/\eps}$ of the semiclassical harmonic oscillator.}
        \label{fig:harmonic_reference_propagator_circuit}
    \end{subfigure}
    \begin{subfigure}[t]{0.95\textwidth}
        \centering
        \[
        \begin{quantikz}[row sep=0.4cm, column sep=0.8cm]
        \lstick{$\ket{m}_{\mathrm d}$}
            & \gate[2]{O_{H_0}(t_m)}
            & \qw
            & \gate[2]{O_{H_0}(-t_m)}
            & \qw
            \\
        \lstick{$\ket{\psi}_s$}
            &
            & \gate[2]{\widetilde O_V}
            &
            & \qw
            \\
        \lstick{$\ket{0}_{\mathrm a}$}
            & \qw
            &
            & \qw
            & \qw
        \end{quantikz}
        \]
        \caption{Circuit for the interaction pictured Hamiltonian 
        $H_{I,d}(t_m)$, where
        $O_{H_0}(\tau)$ denotes the controlled harmonic-oscillator reference evolution $U_{0,d}(\tau)$, and $\widetilde O_V$ is a block-encoding of $V_d$.}
        \label{fig:harmonic_interaction_picture_block}
    \end{subfigure}

    \caption{Quantum circuits for the harmonic-oscillator interaction picture. Top: the circuit for the reference propagator of the semiclassical harmonic oscillator. Bottom: the corresponding HAM-T block obtained by conjugating the potential block with the reference evolution.}
    \label{fig:harmonic_qho_circuits}
\end{figure}
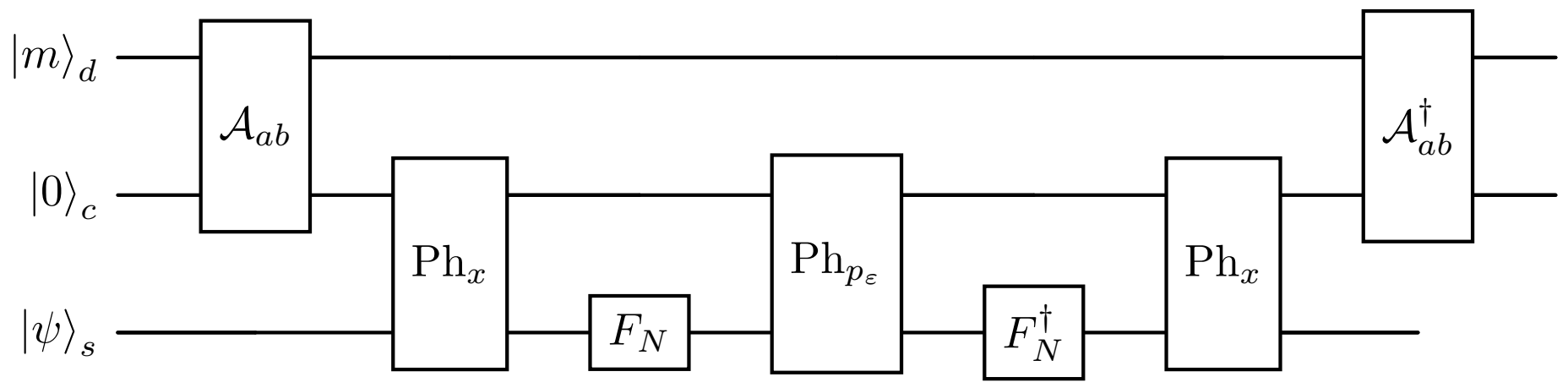
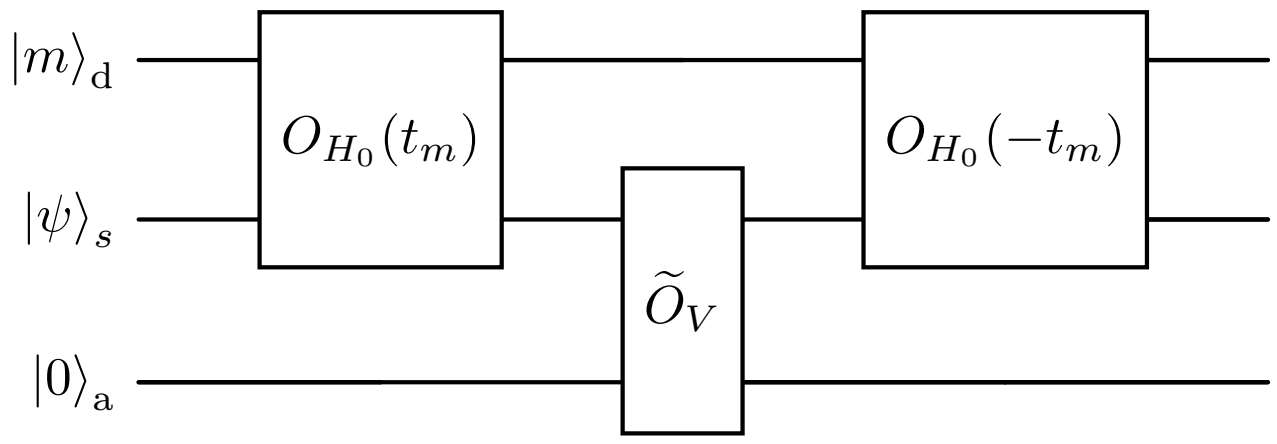

After spatial discretization on a binary grid of size $N=2^n$, this yields the reference unitary $U_{0,d}(t)$ approximating $e^{-itH_0/\eps}$. We use a Fourier spectral discretization. Let $F_N$ denote the discrete Fourier transform from the position basis to the momentum basis. If $x_j$ denotes the position grid and $p_{\eps,k}$ denotes the corresponding semiclassical momentum grid, then one has
\[
U_{0,d}(t)
=
D_x(a(t))\,F_N^\dagger\,D_p(b(t))\,F_N\,D_x(a(t)),
\]
where \[
D_x(a(t))
=
\operatorname{diag}_j\left(e^{-ia(t)x_j^2/\eps}\right),
\qquad
D_p(b(t))
=
\operatorname{diag}_k\left(e^{-ib(t)p_{\eps,k}^2/\eps}\right).
\]
Equivalently, one applies the $x$-quadratic phase, then $F_N$, then the $p_\eps$-quadratic phase, then $F_N^\dagger$, and finally another $x$-quadratic phase, which gives the circuit~\cref{fig:harmonic_qho_circuits}.
Let $\mathcal A_{ab}$ be an arithmetic oracle that evaluates the scalar functions $a(t_m)$ and $b(t_m)$:
\[
\mathcal A_{ab}:
\ket{m}_{\mathrm d}\ket{0}_c
\mapsto
\ket{m}_{\mathrm d}\ket{a_m,b_m}_c,
\qquad
a_m=a(t_m),\quad b_m=b(t_m).
\]
The controlled phase blocks are defined by
\[
\mathrm{Ph}_x:
\ket{a_m,b_m}_c\ket{j}_s
\mapsto
e^{-i a_m x_j^2/\eps}
\ket{a_m,b_m}_c\ket{j}_s,
\quad 
\mathrm{Ph}_{p_\eps}:
\ket{a_m,b_m}_c\ket{k}_s
\mapsto
e^{-i b_m p_{\eps,k}^2/\eps}
\ket{a_m,b_m}_c\ket{k}_s,
\]
in computational and Fourier bases, respectively.
The quadratic phase gates are standard and can be implemented by reversible arithmetic and phase kickback. For example, to implement $D_x(a_m)$ coherently with a time index $m$, we compute a fixed-point approximation of
\[
\phi^{(x)}_{m,j}
=
\frac{a(t_m)x_j^2}{\eps}
\pmod{2\pi},
\]
apply the phase-gradient map
\[
\ket{\phi^{(x)}_{m,j}}\mapsto e^{-i\phi^{(x)}_{m,j}}\ket{\phi^{(x)}_{m,j}},
\]
and then uncompute the work register. The same construction applies to $D_p(b_m)$ in the Fourier basis, using
$
\phi^{(p)}_{m,k}
=
\frac{b(t_m)p_{\eps,k}^2}{\eps}.
$
Therefore, \cref{fig:harmonic_reference_propagator_circuit} implements the oracle for the reference propagator
\[
\ket{m}_d\ket{\psi}_s
\mapsto
\ket{m}_d\,U_{0,d}(t_m)\ket{\psi}_s,
\]
with the work register returned to $\ket{0}_c$.

Each short-time block contains two quantum Fourier transforms and three quadratic phase gates. For a fixed classical time, the phase coefficients can be precomputed, and expanding the squared binary grid indices gives an implementation using $O(n^2)=O(\log^2 N)$ one- and two-qubit gates when arbitrary single-qubit rotations are allowed, consistent with the gate model in~\cite{JainIyerSommaBaoJordan2025}.

For a coherent query over $M$ time points, the cost also includes the arithmetic oracle $\mathcal A_{ab}$ and its inverse. Fix $\omega$, the spatial domain, and the time range $[0,T]$, and let $\delta$ denote the circuit error tolerance for the reference conjugation. Assume that $t_m$ can be computed reversibly from $m$ with cost polynomial in $\log M$ and the requested bit precision. The required coefficient and phase precision is $O(\log N+\log(1/\eps)+\log(1/\delta))$ bits. Standard reversible arithmetic for the time reduction, trigonometric functions, and quadratic phases, together with rotation synthesis, therefore gives a gate cost
\[
\operatorname{poly}\!\left(\log N,\log M,\log(1/\eps),\log(1/\delta)\right)
\]
for the reference part of the interaction-picture block. This estimate concerns circuit implementation; spatial discretization error is separate.

For arbitrary queried times, let $U_{0,d}(t)$ denote the resulting two-block reference circuit. The interaction-picture Hamiltonian defined in~\cref{eq:ip_Ham_def} takes the form
\[
H_{I,d}(t)
=
U_{0,d}(t)^\dagger\,V_d\,U_{0,d}(t),
\quad 
V_d=\operatorname{diag}_j(V(x_j)).
\]
Its circuit is given in \cref{fig:harmonic_interaction_picture_block} and uses two reference-oracle calls and one call to the potential block-encoding $\widetilde O_V$.
Here $V_d$ is block-encoded using a value oracle for the scalar function $V(x_j)$. Namely, assuming $|V(x_j)|\leq \alpha_V$, we use a block-encoding $\widetilde O_V$ satisfying
\[
(\bra{0}_{\mathrm a}\otimes I_s)\,
\widetilde O_V\,
(\ket{0}_{\mathrm a}\otimes I_s)
=
\frac{V_d}{\alpha_V}.
\]
As the matrix is diagonal, such a block-encoding can be obtained by computing $V(x_j)/\alpha_V$, applying a value-controlled rotation on an ancilla qubit, and then uncomputing the value register.

\section{Numerical Experiments}\label{sec:numerics}
In this section, we present two sets of numerical tests. We first verify the nested commutator estimate in
\cref{thm:nested_comm}, focusing on the superconvergent scaling in the time step $h$ for representative nested commutators. We then verify the algorithmic error behavior predicted by \cref{thm:main}, namely the doubled-order convergence in $h$ and the order-independent semiclassical prefactor $\eps^{-1}$. For the algorithmic error, we simulate the quantum Magnus algorithm classically for benchmarking purposes only. We note that such a classical implementation of the quantum circuit is computationally inefficient, with cost scaling as $\mathrm{poly}(N)$ classically, in contrast to its quantum implementation cost as $\mathrm{polylog}(N)$ in the number of spatial degrees of freedom $N$.

We work in one spatial dimension and discretize the operator
\[
A^\eps = -\frac{\eps}{2}\Delta + \frac{1}{\eps}x^2,
\qquad
B^\eps = \frac{1}{\eps}\cos(x),
\]
on a uniform periodic grid with $N$ points on a bounded interval. The Laplacian is approximated by the standard second-order centered finite-difference or spectral discretization operator with periodic boundary conditions. Denoting the resulting matrices by $A_N^\eps$ and $B_N^\eps$, we define the discrete Heisenberg evolution
\[
B_N^\eps(t) := e^{itA_N^\eps} B_N^\eps e^{-itA_N^\eps}.
\]
Since the oscillatory scale of the underlying dynamics becomes finer as $\eps$ decreases, we simultaneously refine the spatial discretization and take $N\propto \eps^{-1}$ so that the discrete model continues to resolve the semiclassical regime. Such a spatial grid refinement is typical and justified for such a semiclassical simulation (see, e.g, \cite{BaoJinMarkowich2002,JinMarkowichSparber2011,LasserLubich2020} for more details).

\subsection{Commutator Scaling in the Time Step and the Semiclassical Parameter}

We first study the dependence on the step size $h$. For this experiment we fix $\eps$ and compute the numerical discretization of the left-hand side of \cref{thm:nested_comm} for a range of step sizes $h$ and various spatial resolutions $N$ large enough. We compare the numerical results with the theoretical bounds.

\begin{equation}\label{eq:numerics-target}
\sup_{\tau_1,\dots,\tau_q\in[t_j,t_{j+1}]}
\max_{C\in\mathscr C_q}
\norm{C(\HI(\tau_1),\dots,\HI(\tau_q))}_{L^2\to L^2}
\le C_{q,T}\eps^{-1} h^{q-1},
\end{equation}
where $\HI(\tau)$ is the spatial discretization of $H_\mathrm{I} (\tau) / \eps$ using a finite-difference matrix discretization on the interval $[-2\pi,2\pi]$, with the operator norm computed as the matrix spectral norm. We verify the superconvergent scaling in the step size $h$, i.e. the factor $h^{q-1}$ for grade $q$ (i.e. $(q-1)$-layer) nested commutators. 

\begin{figure}[tp]
    \centering
    \begin{subfigure}[t]{0.4\textwidth}
        \centering
        \includegraphics[width=\textwidth]{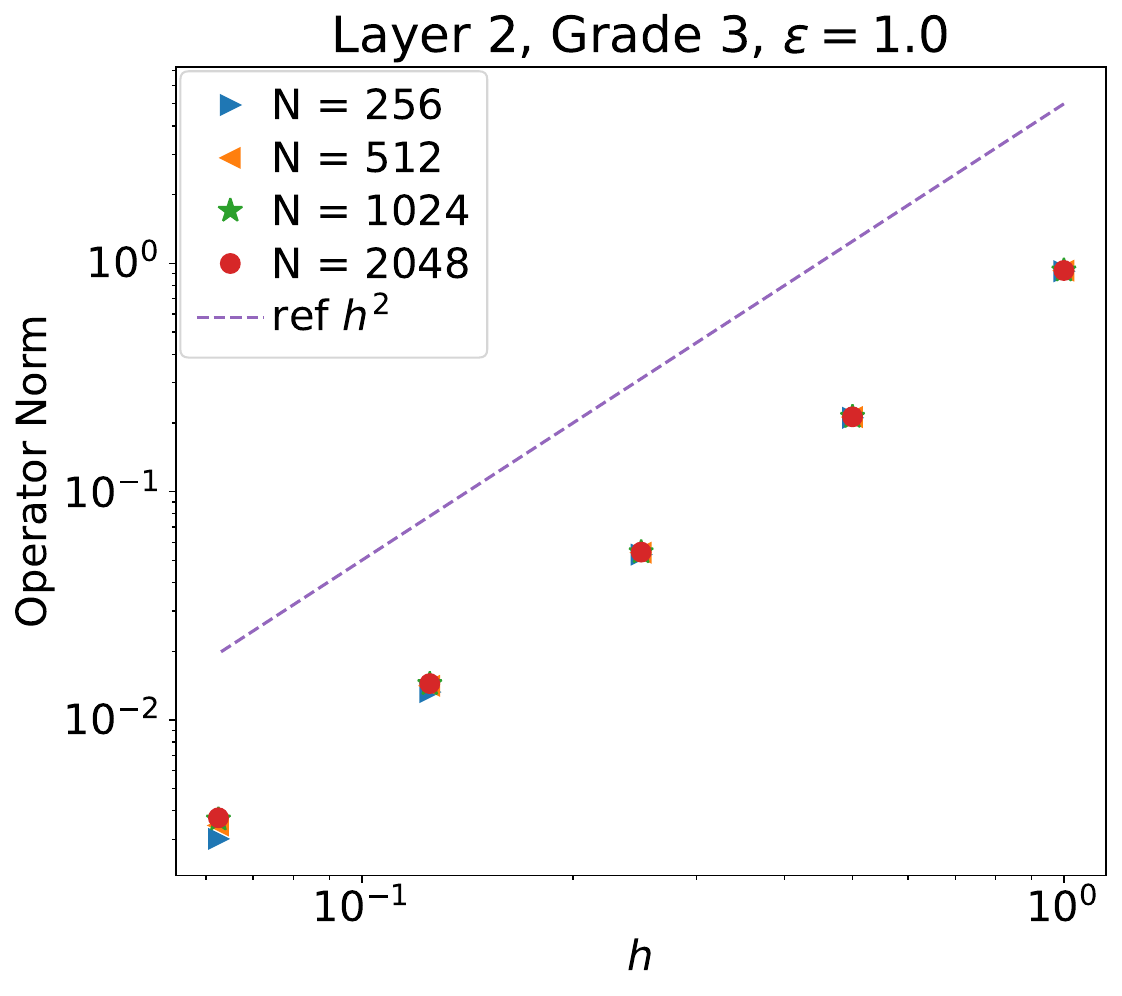}
        \caption{Grade 3 nested commutator.}
        \label{fig:layer2_eps1}
    \end{subfigure}
    \begin{subfigure}[t]{0.4\textwidth}
        \centering
        \includegraphics[width=\textwidth]{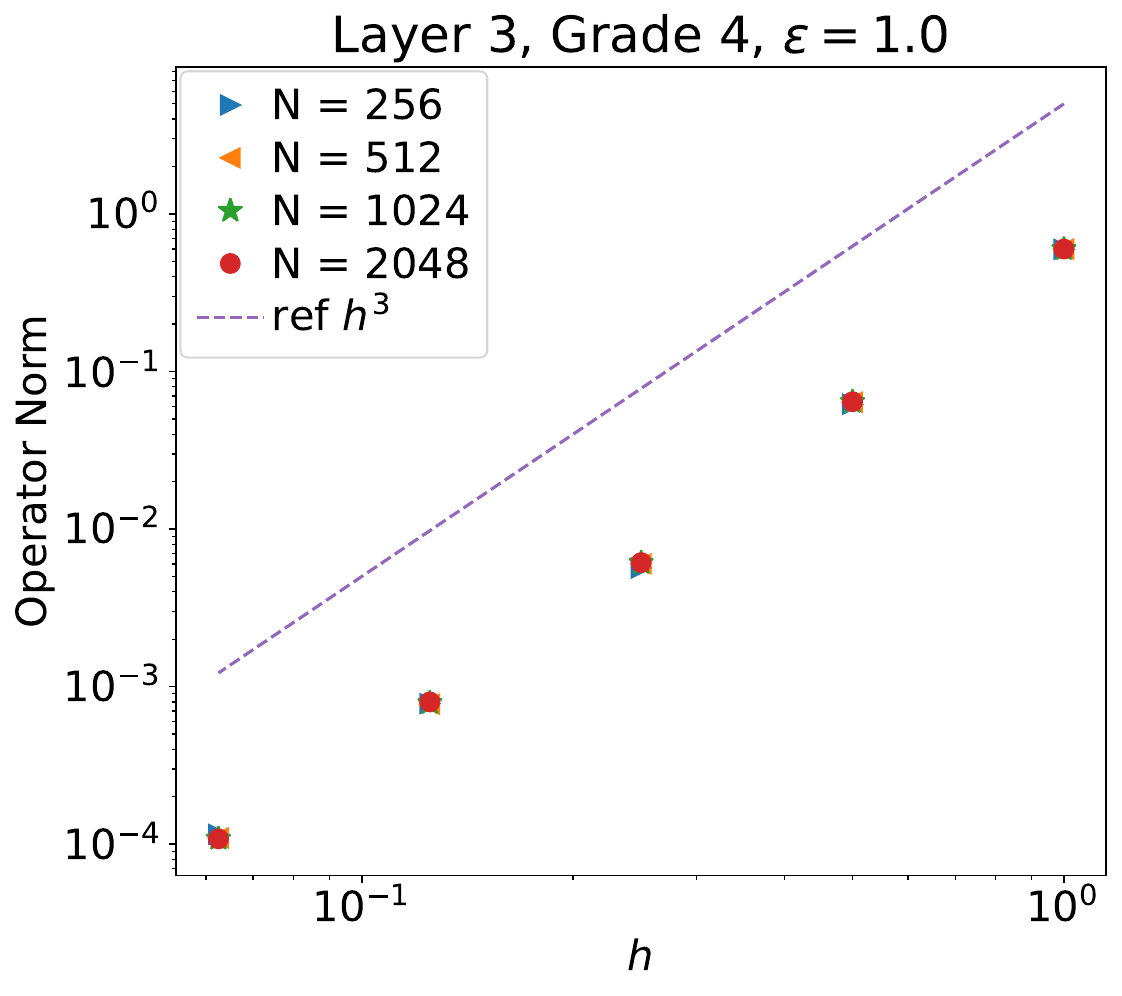}
        \caption{Grade 4 nested commutator.}
        \label{fig:layer3_eps1}
    \end{subfigure}

    \begin{subfigure}[t]{0.4\textwidth}
        \centering
        \includegraphics[width=\textwidth]{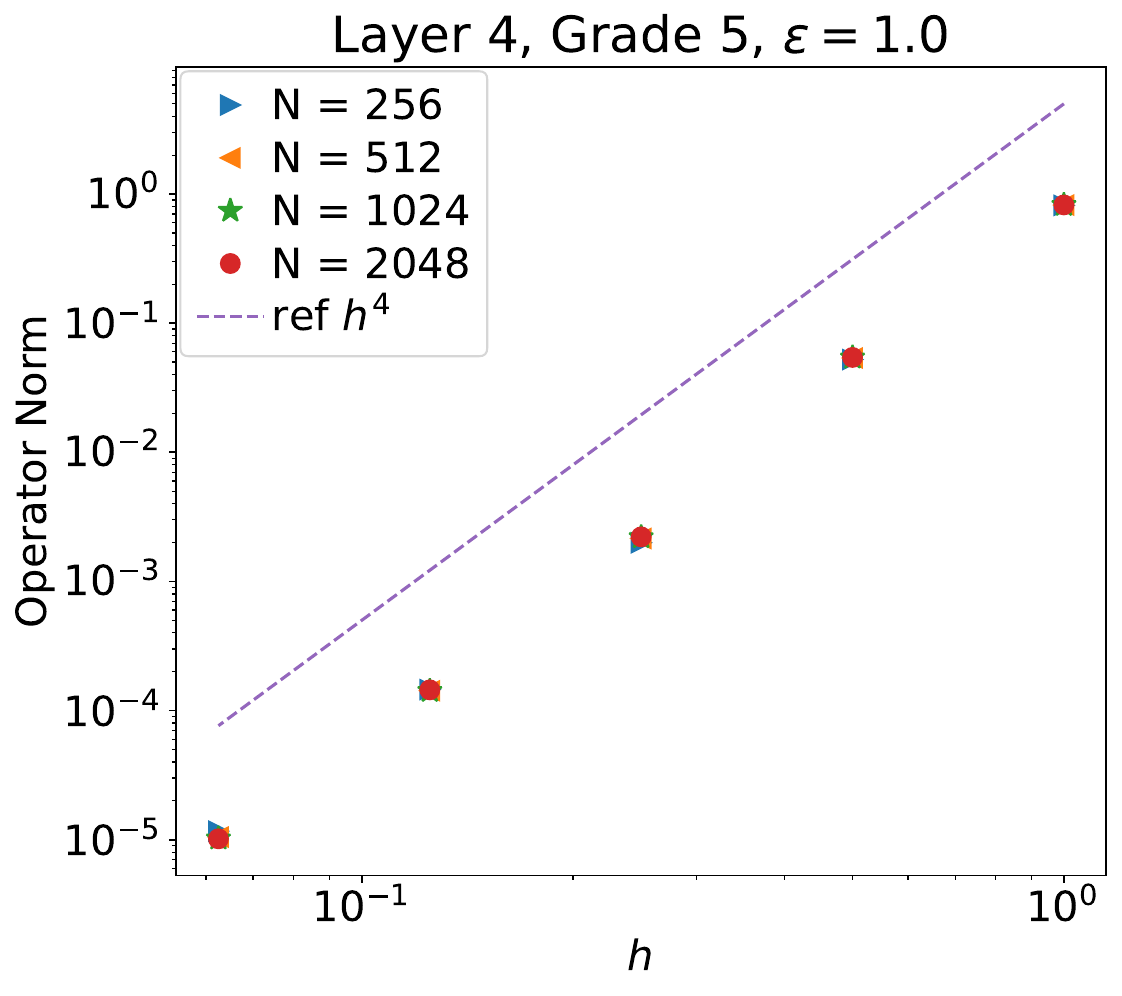}
        \caption{Grade 5 nested commutator.}
        \label{fig:layer4_eps1}
    \end{subfigure}
    \begin{subfigure}[t]{0.4\textwidth}
        \centering
        \includegraphics[width=\textwidth]{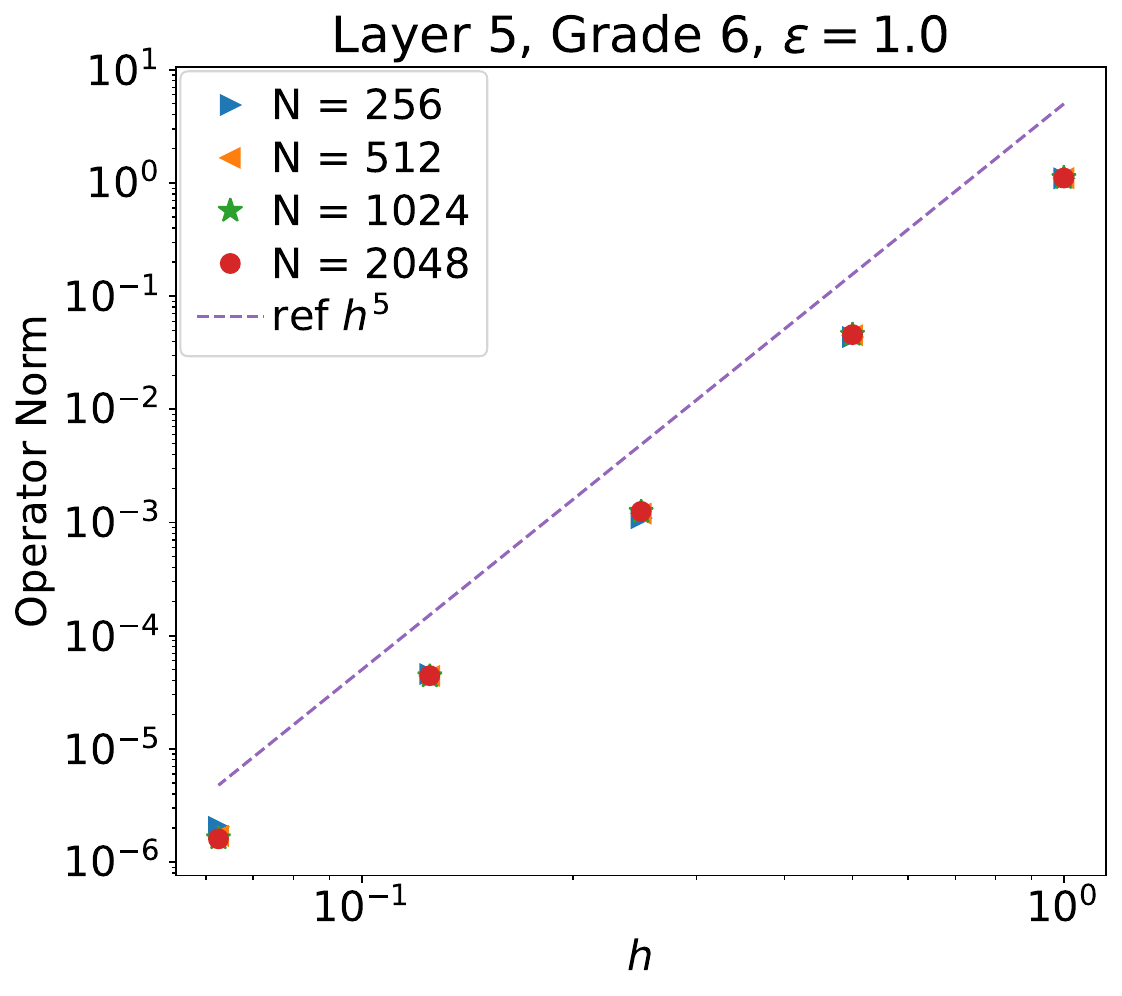}
        \caption{Grade 6 nested commutator.}
        \label{fig:layer5_eps1}
    \end{subfigure}

    \caption{Scaling in the step size $h$ for the grade-$p$ nested commutators with $\eps=1$. The observed slopes agree with the superconvergent scaling predicted by \cref{thm:nested_comm}. The $p$-th grade commutator has an operator norm upper bounded by the $s^{p-1}$ scaling with an error prefactor independent of the spatial grid number $N$.}
    \label{fig:nested_h_eps1}
\end{figure}

We perform this test for the grade $q=3,4,5,6$, corresponding to the commutator depths relevant for the superconvergent local error structure.
We plot the operator norm of the computed representative commutator norms for various time steps
$
h = 2^0,\; 2^{-1},\; 2^{-2},\; 2^{-3},\; 2^{-4},\; $
and various spatial grids 
$
N = 256, 512, 1024, 2048.
$
The log-log plots for $\eps=1$ show that the operator norm scales as $h^{q-1}$, where $q$ is the grade.
\cref{fig:nested_h_eps1} is
in excellent agreement with \cref{thm:nested_comm}. In particular, the observed operator norms are upper bounded by the reference line with a slope that scales as $2, 3, 4, 5$ for the grades $q=3, 4, 5, 6$, respectively. This confirms the superconvergent gain of one power of $h$ per additional commutator layer. Moreover, for sufficiently large $N$, the curves for different spatial resolutions essentially overlap, indicating that the preconstant is not sensitively dependent on the spatial grid $N$, which is important for efficient quantum algorithm implementation with $\mathrm{polylog}(N)$ cost in total. 
This setup strictly extends the previous superconvergence result without harmonic oscillators~\cite{FangZhang2025}.

These computations provide direct numerical evidence for the factor $h^{q-1}$ in \cref{eq:numerics-target}. In the context of the Magnus analysis, this is precisely the mechanism underlying the superconvergent local defect estimates: higher-grade nested commutators gain higher powers of the step size, with the numerics matching the predicted order very closely.

\begin{figure}[tp]
    \centering
    \begin{subfigure}[t]{0.4\textwidth}
        \centering
        \includegraphics[width=\textwidth]{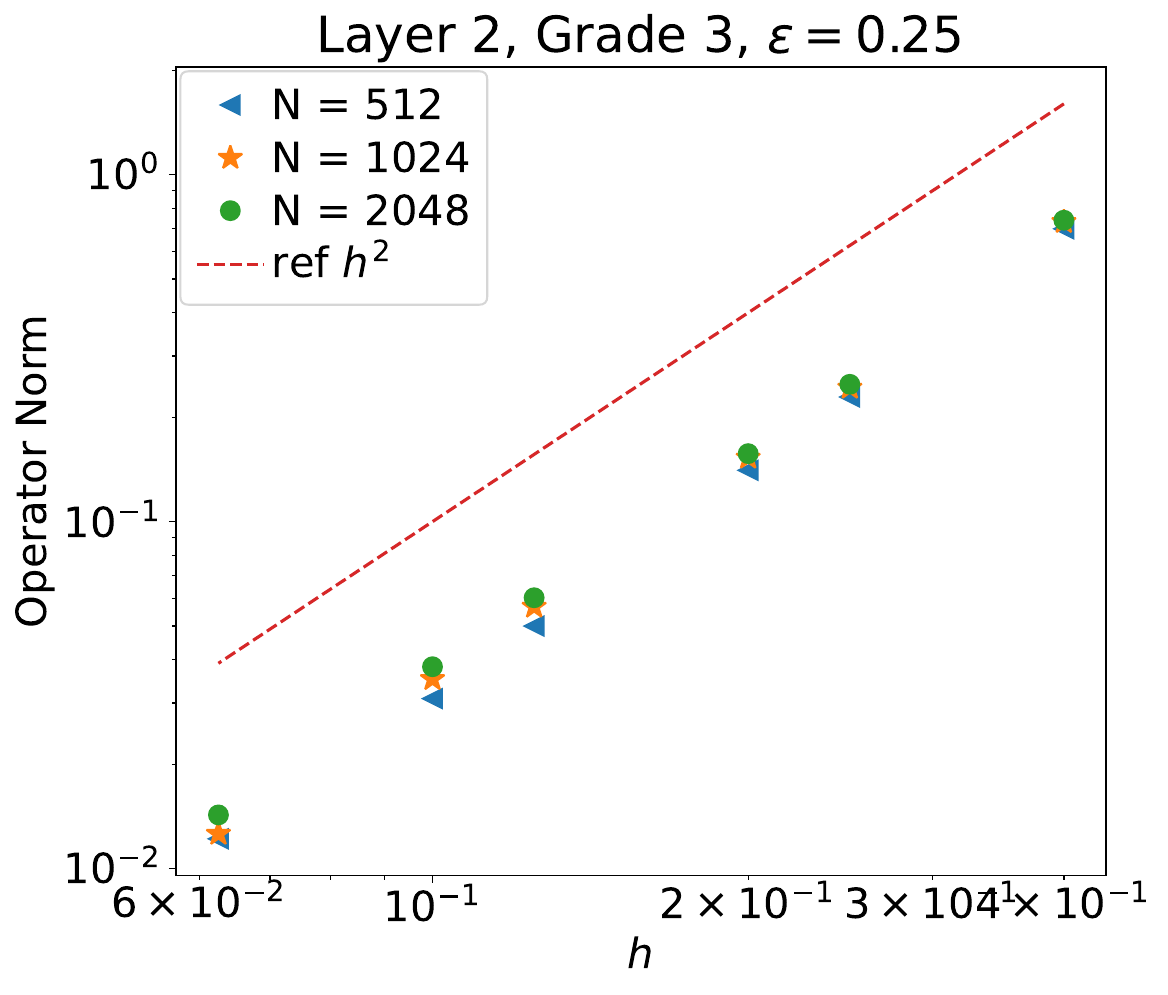}
        \caption{Grade 3 nested commutator.}
        \label{fig:layer2_eps1_4}
    \end{subfigure}
    \begin{subfigure}[t]{0.4\textwidth}
        \centering
        \includegraphics[width=\textwidth]{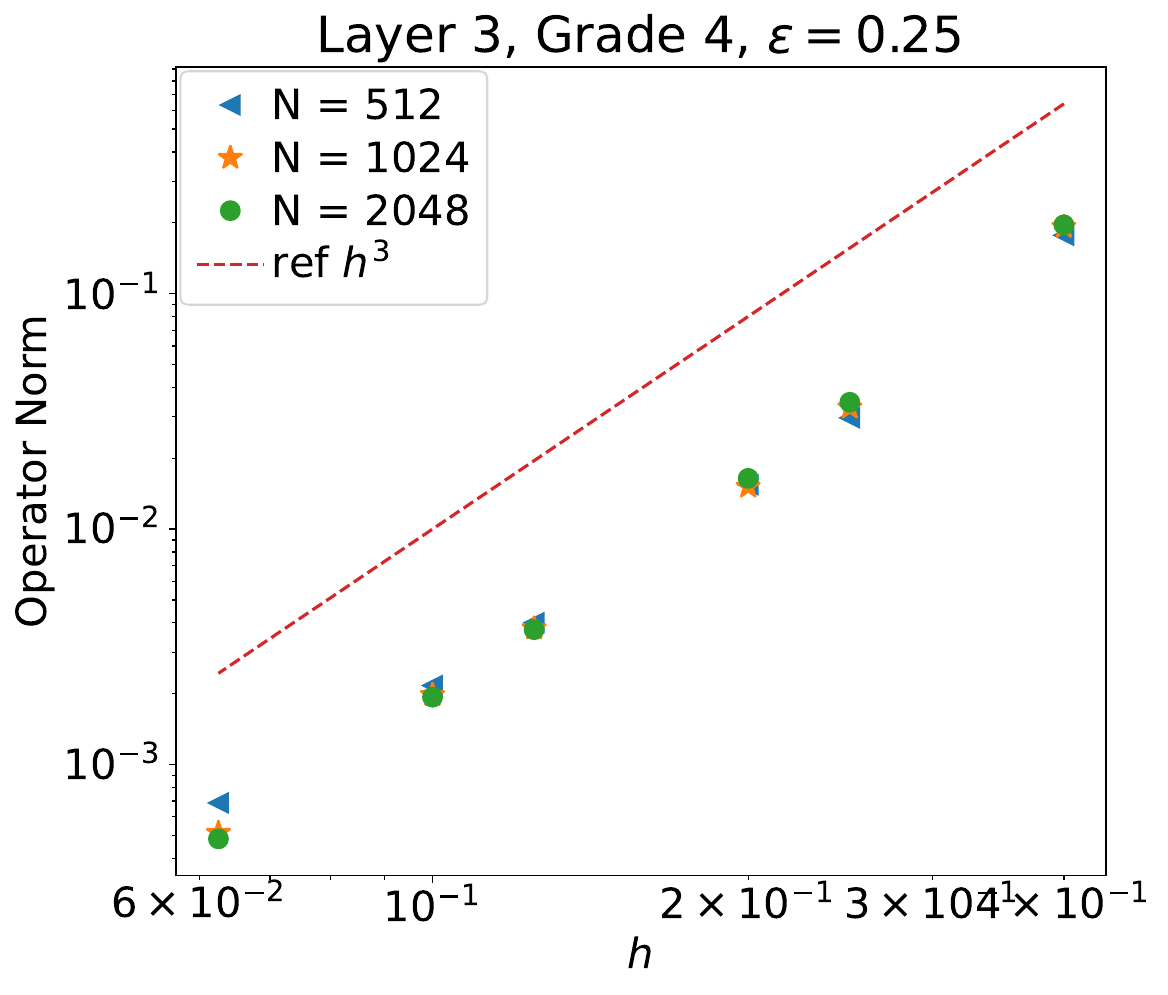}
        \caption{Grade 4 nested commutator.}
        \label{fig:layer3_eps1_4}
    \end{subfigure}

    \begin{subfigure}[t]{0.4\textwidth}
        \centering
        \includegraphics[width=\textwidth]{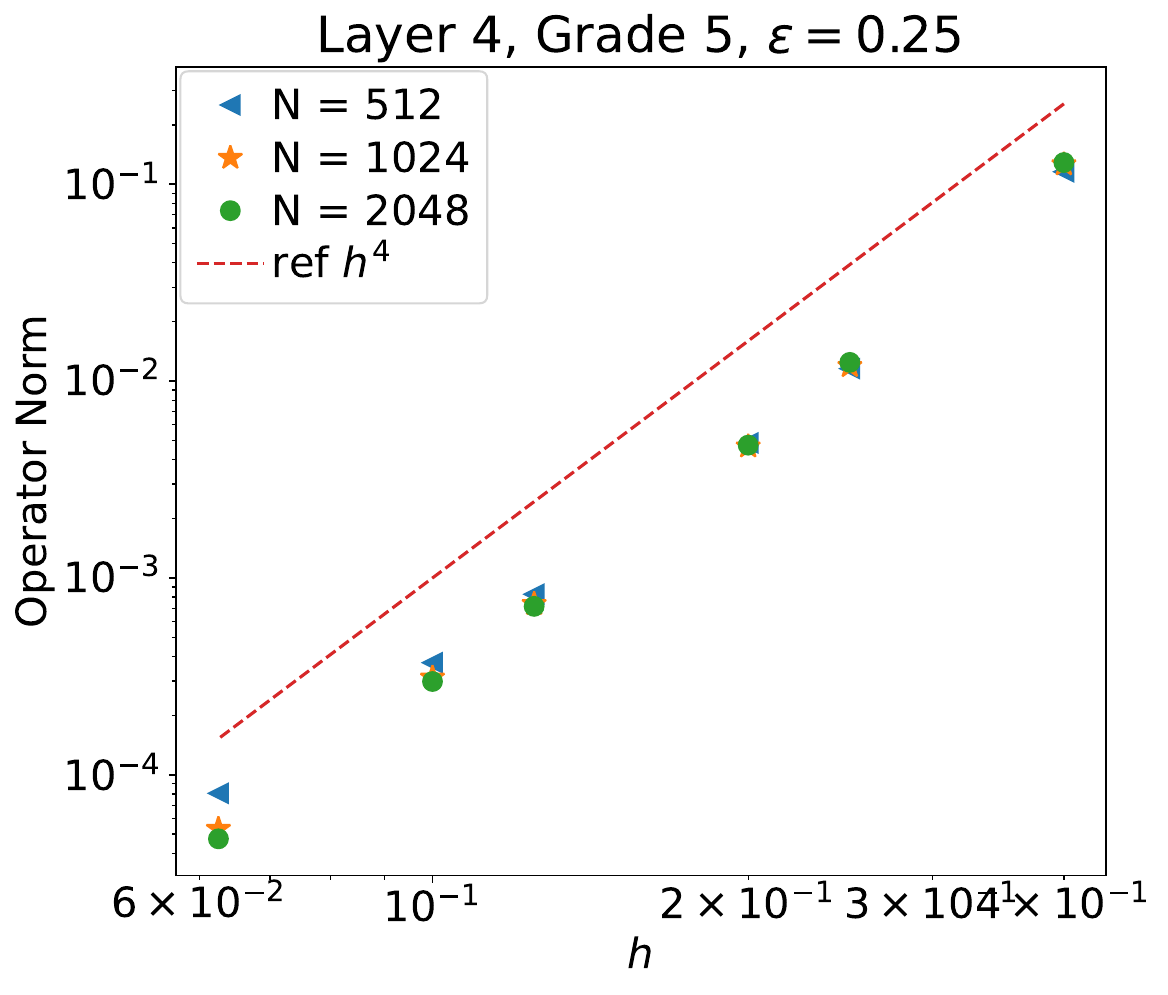}
        \caption{Grade 5 nested commutator.}
        \label{fig:layer4_eps1_4}
    \end{subfigure}
    \begin{subfigure}[t]{0.4\textwidth}
        \centering
        \includegraphics[width=\textwidth]{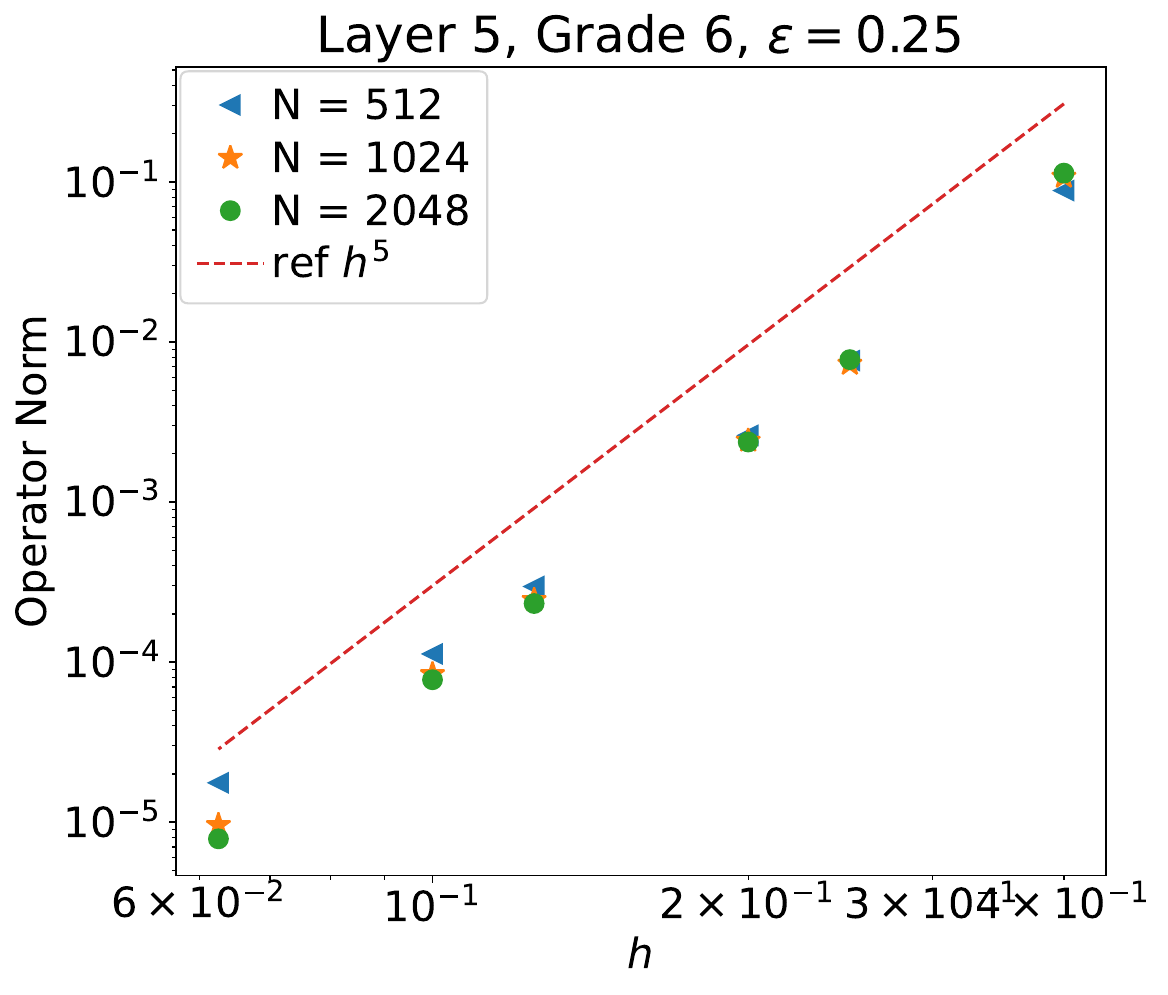}
        \caption{Grade 6 nested commutator.}
        \label{fig:layer5_eps1_4}
    \end{subfigure}

    \caption{Scaling in the step size $h$ for the grade-$p$ nested commutators with $\eps=0.25$. The observed slopes agree with the superconvergent scaling predicted by \cref{thm:nested_comm}. The $p$-th grade commutator has an operator norm upper bounded by the $s^{p-1}$ scaling with an error prefactor insensitive to the spatial grid number $N$ for $N$ large enough resolving the oscillatory dynamics.}
    \label{fig:nested_h_eps1_4}
\end{figure}

\begin{figure}[t]
    \centering
    \begin{subfigure}[t]{0.48\textwidth}
        \centering
        \includegraphics[width=\textwidth]{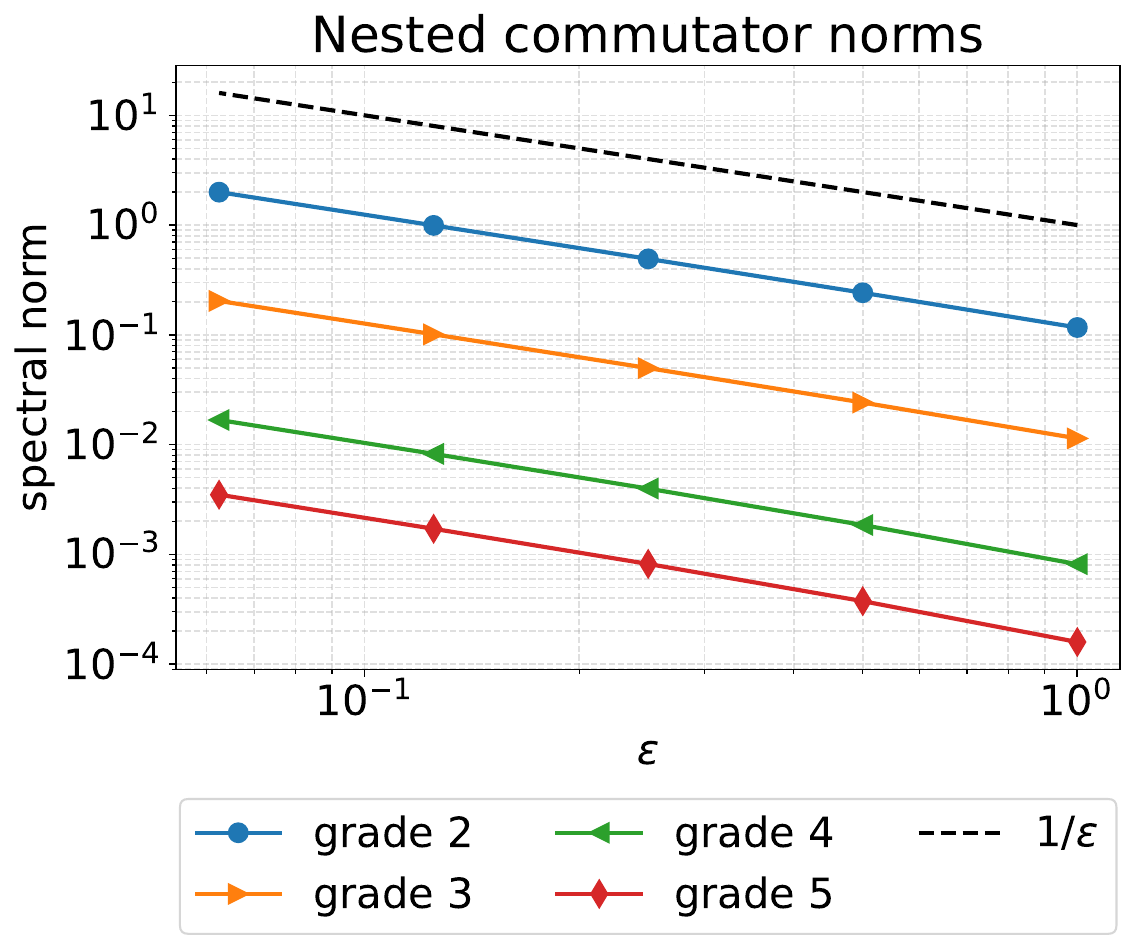}
        \caption{The step size $h = 1/8$.}
        \label{fig:comm_harmonic_eps_dt1_8}
    \end{subfigure}
    \hfill
    \begin{subfigure}[t]{0.48\textwidth}
        \centering
        \includegraphics[width=\textwidth]{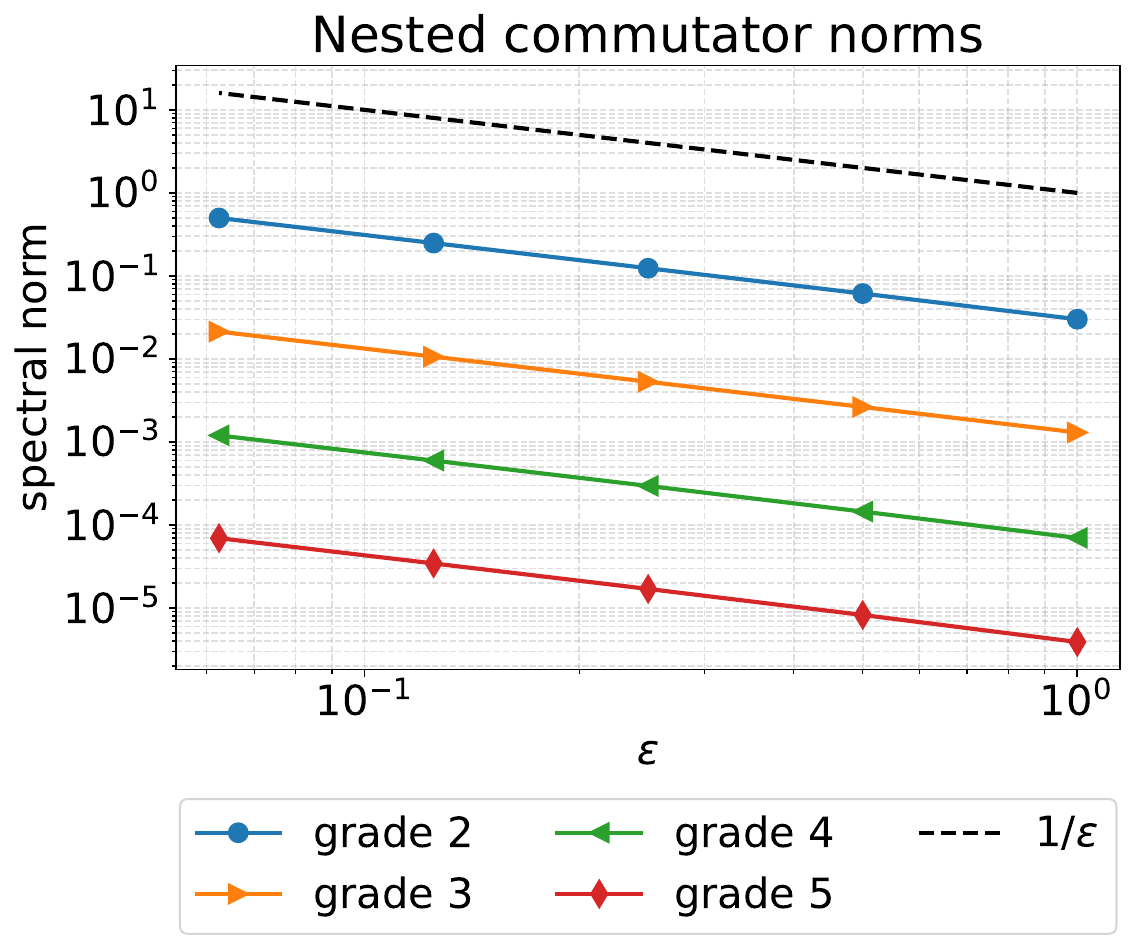}
        \caption{The step size $h = 1/32$.}
        \label{fig:comm_harmonic_eps_dt1_32}
    \end{subfigure}
    \caption{Scaling in the semiclassical parameter $\eps$ for the nested commutators under harmonic confinement, shown for two different time-step sizes. The observed $\eps$-dependence remains of order $\eps^{-1}$ and is insensitive to the commutator grade.}
    \label{fig:comm_harmonic_eps}
\end{figure}

We next study the dependence on the semiclassical parameter $\eps$. For each fixed grade $q$, the log-log plots in \cref{fig:comm_harmonic_eps}
show a clear $\eps^{-1}$ scaling. 
We choose the time-step size $h = 1/8, 1/32$ and plot the representative commutator norms with various $\eps$ values
\[
\eps =  1, 1/2, 1/4, 1/8, 1/16.
\]
Most notably, this behavior is observed uniformly across the tested layers: the $\eps$-dependence does not worsen as the grade $q$ (or the layer) increases. All nested commutator scales are consistent with the single factor $\eps^{-1}$ upper bound predicted by \cref{thm:nested_comm}. 

\subsection{Algorithmic Error for Superconvergence}\label{sec:subsec_num_alg_h}
In what follows, we conduct a classical simulation of the quantum Magnus algorithm. This classical emulator is employed only for benchmarking purposes and is computationally inefficient: the quantum Magnus circuit has cost $\mathrm{polylog}(N)$ in the number of spatial degrees of freedom $N$ if implemented quantumly, whereas its classical simulation requires $\mathrm{poly}(N)$ cost.

\begin{figure}[tp]
    \centering
    \includegraphics[width=0.5\textwidth]{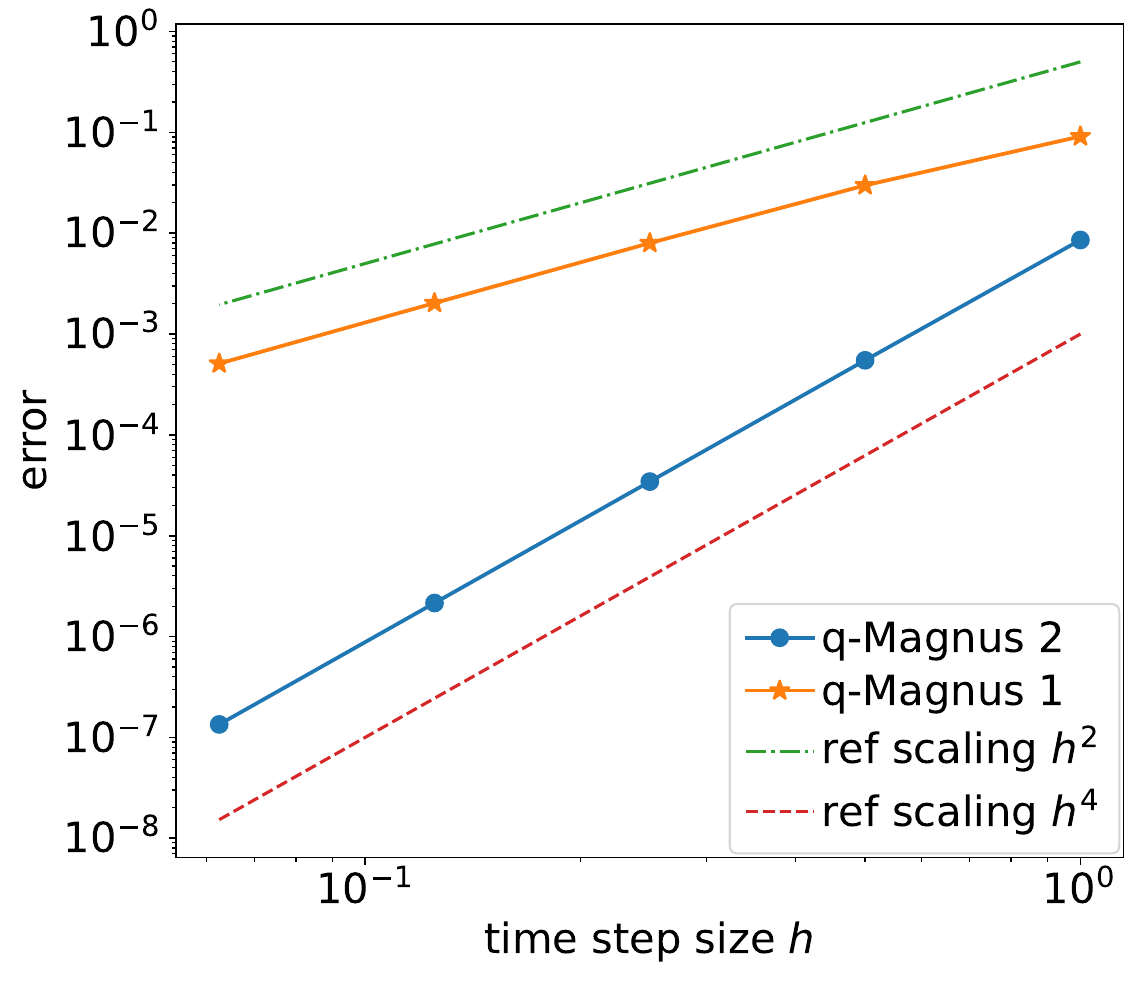}
    \caption{Solution error versus the time-step size $h$ for the first- and second-order quantum Magnus algorithms with harmonic confinement at $T=1$ and $\eps=1/4$. The first-order method exhibits second-order convergence, while the second-order method exhibits fourth-order convergence, in agreement with the superconvergence prediction of \cref{thm:main}.}
    \label{fig:mag12_h_convergence}
\end{figure}

We next verify the solution error behavior predicted by \cref{thm:main}. We first fix $T=1$ and $\eps=1/4$, and choose the number of spatial grid points according to $N=1024$. The solution error is then plotted against the time-step size $h=T/L$, where $L$ ranges over the chosen numbers of time steps. The resulting log--log plot is shown in \cref{fig:mag12_h_convergence}. For both the first-order and second-order quantum Magnus algorithm we use $L=1,2,4,8,16$ and a sufficiently fine quadrature with $n_{\mathrm{quad}}=8192$. The observed slopes clearly show that the first-order Magnus algorithm converges at second order, whereas the second-order Magnus algorithm converges at fourth order. This directly confirms the doubled-order superconvergence predicted by \cref{thm:main} in the harmonic-confinement setting.

\subsection{Algorithmic Error for Semiclassical scaling}\label{sec:subsec_num_alg_eps}

\begin{figure}[tp]
    \centering
    \includegraphics[width=0.5\textwidth]{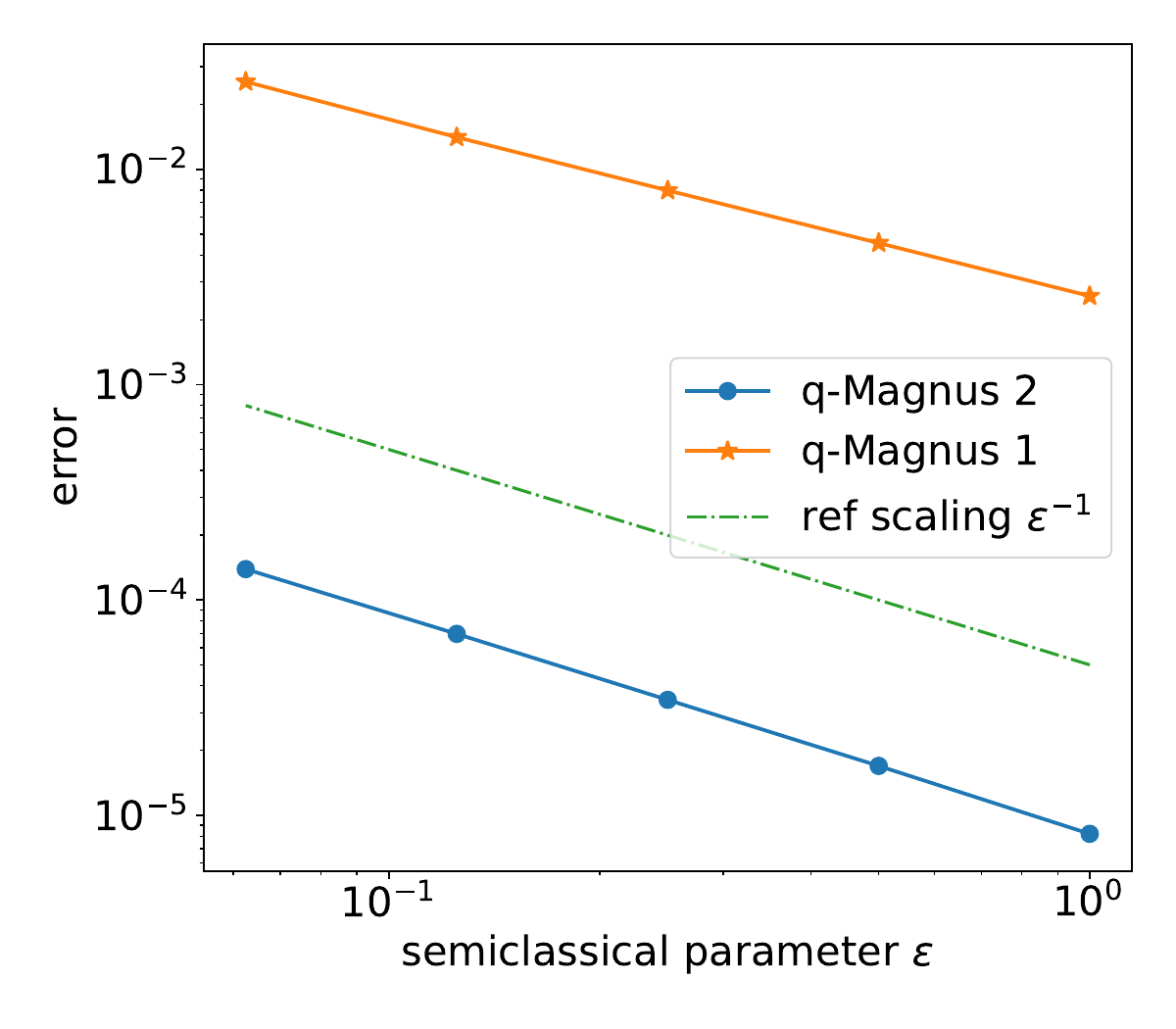}
    \caption{Solution error versus the semiclassical parameter $\eps$ for the first-order and second-order quantum Magnus algorithms with harmonic confinement, fixed $T=1$, and fixed time-step size $h=1/4$. The observed behavior is consistent with the order-independent semiclassical prefactor $\eps^{-1}$ predicted by \cref{thm:main}. For $\eps=1$, the computation is carried out on the enlarged domain $[-3\pi,3\pi)$ in order to capture the less localized wave packet.}
    \label{fig:mag12_eps_scaling}
\end{figure}

We then study the dependence of the algorithmic error on the semiclassical parameter $\eps$. In this experiment we fix $T=1$, the step size $h = 1/4$, and use the first- and second-order Magnus algorithms with $n_{\mathrm{steps}}=4$ and $n_{\mathrm{quad}}=8192$. For smaller values of $\eps$, we again refine the spatial grid according to $N=256/\eps$ in order to resolve the oscillatory dynamics. The initial state is chosen as a normalized semiclassical Gaussian wave packet,
\[
\psi_0(x)
=
C_\eps
\exp\left(-\frac{(x-x_0)^2}{2\eps}\right)
\exp\left(\frac{i p_0(x-x_0)}{\eps}\right),
\qquad
x_0=0,\quad p_0=1,
\]
where $C_\eps$ is the normalization constant ensuring that the $L^2$ norm of the wave function is 1, so that it is a valid quantum state.
For $\eps=1$, however, the wave packet is less localized and propagates farther, so we enlarge the computational domain to $[-3\pi,3\pi)$ and use $N=768$. The resulting plot in \cref{fig:mag12_eps_scaling} shows that the solution error follows the predicted order-independent semiclassical prefactor $\eps^{-1}$, thereby confirming the scaling law in \cref{thm:main}.

Taken together, the numerical experiments confirm both principal features of our analysis. The commutator tests verify the $h^{q-1}$ scaling that underlies the superconvergent local defect estimates, while the solution-error tests confirm the doubled-order global convergence and the order-independent semiclassical prefactor $\eps^{-1}$.

\section{Conclusion and Discussion}

Conceptually, our results show that superconvergence of Magnus quantum algorithms is not a peculiarity of free Schr\"odinger propagation. It persists under a fundamental class of confined semiclassical dynamics, and its true semiclassical scaling is governed by geometric and symbolic structure rather than by formal commutator depth alone. In particular, high-order Magnus methods remain viable in the semiclassical regime to a much greater extent than the abstract commutator expansion would suggest.

Our work continues a recent line of research on Magnus quantum algorithms, but along two different axes. Previous studies established superconvergence for Laplacian-based reference dynamics, extended the phenomenon to arbitrary order, and analyzed discretization-independent or fully discrete settings. The present paper instead addresses the simultaneous issues of structured unbounded reference dynamics and sharp semiclassical prefactors. From the perspective of numerical analysis, these are natural next questions: harmonic-oscillator backgrounds are among the most canonical semiclassical Hamiltonians, and understanding how algorithmic error constants depend on $\eps$ is essential for judging the practical value of a high-order method.

\section*{Acknowledgments}
This work is supported by the National Science Foundation CAREER Award DMS-2438074, and the U.S. Department of Energy, Office of Science, Accelerated Research in Quantum Computing Centers, Quantum Utility through Advanced Computational Quantum Algorithms, grant no. DE-SC0025572. The authors developed and proved the mathematical theory independently by May 2026 without any AI involvement. AI is used for language and grammar checking after the paper was fully written.

\bibliographystyle{unsrt}
\bibliography{magnus-p}
\end{document}